\PassOptionsToPackage{english}{babel}
\documentclass[11pt]{article}

\usepackage{amsmath, amsfonts, amssymb, amsthm,  graphicx, mathtools, enumerate, dsfont}

\usepackage[blocks, affil-it]{authblk}

\allowdisplaybreaks

\usepackage[square]{natbib}
\usepackage[CJKbookmarks=true,
            bookmarksnumbered=true,
			bookmarksopen=true,
			colorlinks=true,
			citecolor=red,
			linkcolor=blue,
			anchorcolor=red,
			urlcolor=blue]{hyperref}
\usepackage[usenames]{color}
\newcommand{\red}{\color{red}}

\usepackage[letterpaper, left=1.2truein, right=1.2truein, top = 1.2truein, bottom = 1.2truein]{geometry}
\usepackage[ruled, vlined, lined, commentsnumbered]{algorithm2e}
\usepackage[font=small,labelfont=it]{caption}

\usepackage{bm}
\usepackage{euscript}
\usepackage{prettyref,soul}
\usepackage[ngerman]{babel}
\usepackage{bbm}
\usepackage{multirow}
\usepackage{diagbox}
\usepackage{hyperref}
    \hypersetup{ colorlinks = true, linkcolor = blue, citecolor = blue}

\let\hat\widehat
\let\bar\overline
\let\tilde\widetilde

\numberwithin{equation}{section} \theoremstyle{plain}
\newtheorem{theorem}{Theorem}[section]
\newtheorem{lemma}{Lemma}[section]
\newtheorem{corollary}{Corollary}[section]
\newtheorem{proposition}{Proposition}[section]

\newtheorem{remark}{Remark}[section]

\def\red{\color{black}}

\DeclareMathOperator{\tr}{Tr}
\DeclareMathOperator{\cov}{Cov}
\DeclareMathOperator{\var}{Var}

\begin{document}
\selectlanguage{english}

\newcommand\E{\mathbb{E}}
\newcommand\EN{\EuScript{N}}
\newcommand\lb{\left(}
\newcommand\rb{\right)}
\newcommand\veps{\varepsilon}
\newcommand\norm[1]{\left\lVert#1\right\rVert}
\newcommand\x{\bm{X}}
\newcommand\y{\bm{Y}}
\newcommand\z{\bm{Z}}
\newcommand\T{\mathsf{T}}
\def\H{\mathbf{H}}

\renewcommand\u{\mathbf{u}}
\renewcommand\a{\mathbf{a}}
\renewcommand{\d}[1]{\ensuremath{\operatorname{d}\!{#1}}}

\newcommand\hf[1]{{\color{black} #1}}

\def\cL{\mathcal{L}}
\def\G{\mathbf{G}}
\def\U{\mathbf{U}}
\def\A{\mathbf{A}}
\def\B{\mathbf{B}}
\def\C{\mathbf{C}}
\def\D{\mathbf{D}}
\def\X{\mathbf{X}}
\def\I{\mathbf{I}}
\def\R{\mathbf{R}}
\def\K{\mathbf{K}}
\def\S{\mathbf{S}}
\def\V{\mathbf{V}}
\def\Y{\mathbf{Y}}
\def\diag{\text{diag}}
\def\s{\underline{s}}
\def\m{\underline{m}}

  \title{Asymptotic joint distribution of extreme eigenvalues and trace of large sample covariance matrix in a generalized spiked population model}
 
\author[1]{Zeng Li}
\author[2]{Fang Han}
\author[3]{Jianfeng Yao}
\affil[1]{
Pennsylvania State University
 \\
zxl278@psu.edu
}
\affil[2]{
University of Washington
 \\
fanghan@uw.edu
}
\affil[3]{
The University of Hong Kong
 \\
jeffyao@hku.hk
}

\date{}

\maketitle

\begin{abstract}
This paper studies the joint limiting behavior of extreme eigenvalues and trace of large sample covariance matrix  in a generalized spiked population model, where the asymptotic regime is such that the dimension and sample size grow proportionally. 
The form of the joint limiting distribution is applied to conduct Johnson-Graybill-type tests, a family of approaches testing for signals in a statistical model. For this, higher order correction is further made, helping alleviate the impact of finite-sample bias. The proof rests on determining the joint asymptotic behavior of two classes of spectral processes, corresponding to the extreme and linear spectral statistics respectively. 
\end{abstract}



\section{Introduction}\label{sec:intro}




Considering a sequence of independent and identically distributed (i.i.d.) $p$-dimensional real-valued random vectors $\{\x_1,\ldots,\x_n\}$ with zero mean and population covariance matrix $\mathbf{\Sigma}_p$,
{\red the corresponding sample covariance matrix is defined as} 
\begin{equation}\label{eq:sn}
\S_n=\frac 1n\sum_{i=1}^n \x_i\x_i^{\mathsf{T}},
\end{equation}
with $\lambda_1\geq\lambda_2\geq \cdots\geq \lambda_p\geq0$ denoting the eigenvalues of $\S_n$. It is statistically fundamental and important to study the distributions of the $m$ largest eigenvalues $\lambda_1,\ldots,\lambda_m$ and the trace, $\tr(\S_n)=\sum_{j=1}^p\lambda_j$, of $\S_n$ as $p=p_n$ grows to infinity with $n$. 

Indeed, each of them has led to a large volume of literature. For results on extreme sample eigenvalues, \cite{Johnstone01} first introduced the spiked population model as the non-null case where all eigenvalues of $\mathbf{\Sigma}_p$ are unit except for a fixed small number of spikes, i.e., 
\begin{equation}\label{eq:Johnstone}
\mbox{Spec}(\mathbf{\Sigma}_p)=\Big\{\alpha_1,\cdots,\alpha_m,\underbrace{1,\cdots, 1}_{p-m}\Big\}.
\end{equation}
Here we define $\mbox{Spec}(\A)$ to be the sets of eigenvalues of matrix $\A$.
Under the ``null" (i.e., $\mathbf{\Sigma}_p$ is the identity {\red matrix} $\I_p$), \cite{Johnstone01} established the Tracy-Widom law  for the {\red largest eigenvalue} of real Wishart matrix $\S_n$. Following  Johnstone's development, many efforts have been put into quantifying the effect caused by spiked eigenvalues $\{\alpha_k,~1\leq k\leq m\}$ on $m$ extreme sample ones $\{\lambda_k,~1\leq k\leq m\}$. To name a few, under Johnstone's spiked model settings, \cite{BS06} thoroughly studied the almost sure limits of the extreme sample eigenvalues under the {\em Mar\v{c}enko-Pastur regime} 
when $p,n\rightarrow \infty,~p/n\rightarrow y\in(0,\infty)$. They found that these limits are different when the corresponding population spiked eigenvalues are larger or smaller than critical values $1+\sqrt{y}$ and $1-\sqrt{y}$. Similar phase transition phenomenon of largest sample eigenvalues was shown in \cite{BBP05} for complex Gaussian population. \cite{Paul07} further proved that a phase transition of eigenvectors also occurs with Gaussian observations. \cite{BY08} followed the set-up of \cite{BS06} and established central limit theorems (CLTs) for the extreme sample eigenvalues associated with spikes outside the interval $[1-\sqrt{y},1+\sqrt{y}]$ under general population distributions. \cite{BY12} extended the theory in \cite{BY08} to a generalized spiked population model where the base population covariance matrix is arbitrary.

In contrast to extreme sample eigenvalues, many important statistics in multivariate analysis can be expressed as linear functionals of eigenvalues of some random matrices, namely, linear spectral statistics (LSS). $\tr(\S_n)$ is one of the most important examples. Limiting behaviors of LSS has been intensively studied in the literature. 
One of the most widely used results is \cite{BS04}, which first established the asymptotic normality for LSS of sample covariance matrix $\S_n$ under the  Mar\v{c}enko-Pastur regime with some moment restrictions on data entries. Further refinement and extensions can be found in numerous follow-up works. To name a few, \cite{Zheng15} studied {\red  CLT for LSS of} sample covariance matrix when the population mean vector is unknown. \cite{CP15} investigated the ultra-high dimensional case when the dimension $p$ is much larger than the sample size $n$. They further established the asymptotic normality for LSS as $p/n\rightarrow \infty$ and $n\rightarrow \infty$. \cite{ZBYZ16}  removed the fourth order moment condition in \cite{BS04} and incorporated it into the limiting parameters. \cite{ZBY-cltF} derived a  CLT for LSS of large dimensional general Fisher matrices. The limiting distribution of $\tr(\S_n)$ is derivable by implementing these results.

\hf{Despite the substantial advances in both directions, to our knowledge, little has been made on investigating the joint distribution of extreme sample eigenvalues and trace, which is equivalent to studying the asymptotic joint {\red distribution of the largest and summation of sample eigenvalues}. As will be seen soon, obtaining such a limiting distribution is fundamental in many applications, and is worth investigating in depth.} 

\hf{As a first contribution of this paper, we aim to study such a joint distribution. For this, we {\red focus} on Bai and Yao's generalized spiked model \citep{BY12}, which generalizes Johnstone's spiked model {\red in} \cite{Johnstone01}.} Here, the population covariance matrix $\mathbf{\Sigma}_p$ has the structure
\begin{equation}\label{eq:GeneralSpikeModel}
\mathbf{\Sigma}_p=\lb
\begin{array}{cc}
{\bf\Lambda}&{\bf 0}\\
{\bf 0} & {\bf V}_{p'}
\end{array}\rb,
\end{equation}
where ${\bf \Lambda}$ and ${\bf V}_{p'}$ are of dimension $m\times m$ and $ p'\times p'~(p'=p-m)$ and ${\bf \Lambda}$ is assumed fixed.  The eigenvalues of ${\bf \Lambda}$ and ${\bf V}_{p'}$ are \hf{$\alpha_1\geq\cdots\geq\alpha_m>0$} and $\beta_{p',1}\geq \cdots\geq \beta_{p',p'}\geq 0$, respectively.  The $\alpha_j's$ are larger than and well separated from ${\beta_{p',j}}'s$, thus named as {\em spiked eigenvalues}.

{\red Under the generalized spiked model with $p/n\to y\in(0,\infty)$, we prove that the extreme eigenvalues and $\tr(\S_n)$ are jointly asymptotically normal and asymptotically independent.} The results are hence connected to the influential work of \cite{CT78} and \cite{Hsing95} on sum and maximum of i.i.d. and strongly mixing random variables. The conclusion holds as long as finite fourth order moments {\red exist}, and in particular, requires no normality assumption.   Our result is hence also connected to {\red another} interesting related work \citep{Davis16} where for sample covariance matrix $\S_n$ with heavy-tailed entries, this asymptotic independence also holds.

Although $\tr(\S_n)$ can be represented {\red as} the summation of all sample eigenvalues, in fact it is very difficult to quantify the correlation between extreme eigenvalues and the rest {\red bulky} ones, especially under the high dimensional settings without a Gaussian assumption. In facing this challenge, we make full use of the spiked model structure and carry out a block-decomposition analysis of spiked and non-spiked ones. The correlation between extreme eigenvalues and trace of each block of $\S_n$ is analyzed separately based on the joint asymptotic behavior of two classes of spectral processes, corresponding to the extreme and linear spectral statistics respectively. The idea of block-decomposition  provides a novel perspective for proving the asymptotic independency between $\lb\lambda_k\rb_{1\leq k\leq m}$ and $\tr(\S_n)$.  \cite{Paul07} adopted a similar block-decomposition technique to represent sample eigenvalues as solutions to certain equations. However, \cite{Paul07} only considered the Gaussian case and focused on individual behavior of eigenvalues and eigenvectors when their population covariance is spiked with unit bulk eigenvalues.

\hf{The form of joint limiting distribution is then employed to conduct Johnson-Graybill-type tests \citep{JG72}, a family of approaches testing for signals in a statistical model based on the sample ratio $\lambda_1/(\tr(\S_n)/p)$. To name an application, this family of tests is important in} modern signal processing applications, such as testing for the presence of signals in cognitive radio and non-parametric signal detection in array processing. \hf{For more details, we refer the readers to \cite{BDMN11} and a comprehensive review paper \citep{PA14}}.  \cite{Johnstone01} proved that when $\mathbf{\Sigma}_p=\sigma^2\I_p$,  $\lambda_1/\sigma^2$ converges to the Tracy-Widom (TW) distribution after appropriate centering and scaling.
However, this ratio test statistic $\lambda_1/(\tr(\S_n)/p)$ cannot be well approximated by the same asymptotic distribution. Finite sample adjustment of critical values for every $(p,n)$ combinations was made in \cite{Nadler11}. \cite{Deo16} suggested an alternative variance correction which also improved the finite sample approximation. However, neither of them derived the asymptotic distribution of this studentized statistic. Furthermore, in the presence of $m$ spikes, the limiting TW distribution of $\lambda_{m+1}/\sigma^2$ has not been fully testified yet. \cite{Deo16} performed some simulation experiments of testing multiple spikes based on the TW conjecture of $\lambda_{m+1}/\sigma^2$. They found that the test was uniformly undersized due to the downward bias of the TW approximation. Further analytic tools are needed to correct this bias, \hf{which is nevertheless nontrivial}. 

In this paper we start from a different perspective by studying the behavior of $\lambda_m/\sigma^2$ in the presence of $m$ spikes. {\red We use $(p-m)^{-1}\sum_{j>m}\lambda_j$ instead of $p^{-1}\tr(\S_n)$ as the surrogate for $\sigma^2$. Although the analytic tools are the same, the former enjoys better performance in finite sample cases.} As a specific example, we formulate our null hypothesis as
 the spiked covariance model where
\begin{equation}\label{eq:SpikeCov}
\mathbf{\Sigma}_p=\sum_{i=1}^m a_{i}\bm{ v}_i\bm{v}_i^{\T}  +\sigma^2\I_p,
\end{equation}
and $\bm{v}_i$'s are orthonormal vectors.  A typical example of such parametrization originates from the factor model where the $p$-dimensional  data vector $\x_t$ has a factor structure of the form
\[\x_t=\A \bm{F}_t+\bm{E}_t,\]
with $\bm{E}_t\sim \EN(0,\sigma^2\I_p)$ independent of $\bm{F}_t\sim\EN(0,\I_m)$ and $\A_{p\times m}$ is a deterministic matrix such that $\A\A^{\T}  $ has spectrum $\sum_{i=1}^m \alpha_{i}\bm{v}_i\bm{v}_i^{\T}$.
The limiting distribution of our test statistic{\red,  $\lambda_m/\frac{1}{p-m}\sum_{j>m}\lambda_j$,} is then derived based on the asymptotic joint distribution of $\lb\lb\lambda_k\rb_{1\leq k\leq m},\tr(\S_{n})\rb^\hf{\T}$, and {\red the corresponding test is implementable} due to the developed theory. Our test targets at detection of signals above certain signal-to-noise ratio. Higher order corrections are further made to alleviate the impact of finite sample bias, which ensures satisfactory testing size and power. 

It is worth mentioning here that this test is closely related to sphericity test (i.e. to test $H_0:m=0,~\mbox{v.s.}~H_1:m>0$) discussed in \cite{Onatski13} and \cite{WZY14}. In particular, they have non-zero power under the spiked alternative (i.e. $H_1:m>0$) even when the spikes are below the phase transition threshold. Admittedly, our test statistic can only detect distant spikes above the phase transition threshold. However, the sphericity test in \cite{Onatski13} and \cite{WZY14} is only designed for testing existence of signals while our tests can be used to detect total number of spikes and the signal strength of the spikes {\red being} tested. In another related work, \cite{Choi17} formulated the observed data matrix $\X\in \mathbb{R}^{n\times p}$ as the sum of a low-rank signal matrix $\B\in \mathbb{R}^{n\times p}$ and a Gaussian noise matrix $\mathbf{E}\in \mathbb{R}^{n\times p}$ and aimed at finding the rank of the deterministic signal matrix $\B$. It is, however, very different from our model settings. {\red In their model signals are treated as a low-rank mean of the observed data matrix} while in this paper we analyze spiked models when factors are embedded in a spiked population covariance structure.

Throughout the paper, we use bold Roman capital letters to represent matrices, e.g., $\A$. $\tr(\A)$ and $|\A|$ denote the trace and determinant of matrix $\A$.  For matrix $\A$, $[\A]_{ij}$ denotes the $(i,j)$-th entry of $\A$. $\mbox{diag}\lb\alpha_1,\cdots,\alpha_m\rb$ represents an $m\times m$ diagonal matrix with diagonal entries $\alpha_1,\cdots,\alpha_m$. Scalars are often in lowercase letters and random ones in capitals.  Vectors follow bold italic style like $\bm{v}_i$ and random vectors are in capitals like $\bm{F}_i$. $\mathbb{N}$, $\mathbb{R}$, and  $\mathbb{C}$ represent  the sets of natural, real, and complex numbers.  $\bm{1}(\cdot)$ stands for indicator function and $\T$ stands for transpose of vectors or matrices. Let $f:\mathbb{C}\rightarrow \mathbb{C}$ be a complex-valued function defined on the complex plane $\mathbb{C}$, then $\oint_{\gamma} f(z)\d z$ denotes the contour integral of $f(z)$ on the Jordan curve $\gamma$. For any $x\in\mathbb{R}$, $\delta_x$ represents the point mass at $x$.

The remaining sections are organized as follows. Section \ref{sec:Prem} gives a detailed description of the generalized spiked model and introduces some preliminary results which form the basis of our analysis. Our main results are presented in Section \ref{sec:Joint}. An application to factor modeling is studied in depth in Section \ref{sec:App}. Proofs of theorems and technical lemmas are relegated to Section \ref{sec:proofs}. 

\section{Generalized spiked population model and preliminaries}\label{sec:Prem}

For any $p\times p$ square matrix $\A$ with eigenvalues $(\theta_j)_{1\leq j\leq p}$, its empirical spectral distribution (ESD) is the measure $F^\A=p^{-1}\sum_{j=1}^p\delta_{\theta_j}$ (weighting equally the eigenvalues). Under the generalized spiked population model \eqref{eq:GeneralSpikeModel}, \hf{the following assumptions are made:}
\begin{itemize}
	\item[(i)] as  $n\rightarrow \infty$, $p=p_n\rightarrow \infty$ such that $p/n\rightarrow y\in(0,\infty)$;

	\item[(ii)] the sequence of spectral norms of $\mathbf{\Sigma}_p$ is bounded and the ESD $H_{p'}$ of $\V_{p'}$ converges to a nonrandom limiting distribution $H$;
	
	\item[(iii)] the eigenvalues $\{\beta_{p',j},~1\leq j\leq p'\}$ of $\V_{p'}$ are such that as $n\rightarrow \infty$,
	\[\sup_{j\leq p'} d(\beta_{p',j},\Gamma_H)=\veps_p\rightarrow 0,\]
	where $d(x,A)$ denotes the Euclidean distance of $x$ to a set $A$ and $\Gamma_H$ \hf{stands for} the support of $H$;
	
	\item[(iv)] the sample vectors $\x_i,~1\leq i\leq n$ can be expressed as $\x_i=\mathbf{\Sigma}_p^{1/2}\y_i$, where $\y_i$ are i.i.d. $p$-dimensional vectors with i.i.d. components $\{Y_{ij},~j=1,\cdots,p\}$ satisfying $\E Y_{ij}=0$, $\E |Y_{ij}|^2=1$, and $\E |Y_{ij}|^4=\nu_4<\infty$.
	
\end{itemize}

\noindent
Letting $\mu$ be a finite measure on the real line with support $\Gamma_\mu$, its Stieltjes transform $s_\mu(z)$ is defined as
\[s_\mu(z)=\int \frac{1}{x-z}\mu (\d x),~z\in \mathbb{C}^+,\]
where $\mathbb{C}^+:=\{z\in\mathbb{C}:\Im(z)>0\}$ is the upper half plane with positive imaginary part \hf{and $\Im(z)$ denotes the imaginary part of any given complex value $z$}. 

Let  $F^{\S_n}$ be the ESD of the sample covariance matrix $\S_n$. It is well known that under Assumptions (i) to (iv), $F^{\S_n}$ weakly converges to a nonrandom probability measure $F^{y,H}$, the Mar\v{c}enko-Pastur (M.P.) distribution with indexes $(y, H)$. Its Stieltjes transform $s(z)$ is implicitly defined as a solution to the equation
\[s(z)=\int \cfrac{1}{t(1-y-yzs(z))-z}\d H(t).\]
Correspondingly, the Stieltjes transform $s_n(z)=\frac 1p\tr\lb \S_n-z{\bf I}_p\rb^{-1}$ of the ESD $F^{\S_n}$ converges to $s(z)$ almost surely as $n\rightarrow \infty$.

Moreover, consider  \hf{an} $n\times n$ companion matrix for $\S_n$,
\[\underline{\S}_n:=\frac1n \Y^{\T}  \mathbf{\Sigma}_p \Y~~~{\rm with}~~~  \Y=(\y_1,\cdots,\y_n).\]
Both matrices share the same non-null eigenvalues and their ESDs satisfy 
\[nF^{\underline{\S}_n}-pF^{\S_n}=(n-p)\delta_0.\]
Their limits and respective Stieltjes transforms are linked to each other by the relation
\[\underline{F}^{y,H}-yF^{y,H}=(1-y)\delta_0,\quad\underline{s}(z)=-\frac{1-y}{z}+ys(z),\]
and the Stieltjes transform $\underline{s}(z)$ of $F^{\underline{\S}_n}$ satisfies the Silverstein equation \citep{SilvChoi95}:
\begin{equation}\label{eq:st}
z=-\frac{1}{\s}+y\int\frac{t}{1+t\s}\d H(t).
\end{equation}

Notice that the spiked structure \eqref{eq:GeneralSpikeModel} can be viewed as a finite rank perturbation of a general population covariance matrix with eigenvalues $\{\beta_{p',j}\}$. As the number of spikes $m$ is fixed while $p\rightarrow \infty$, the limiting spectral distribution of $F^{\S_n}$ is determined by the distribution of bulk population eigenvalues $\{\beta_{p',j}\}$ independent of the spikes. However, the behavior of the $m$ extreme sample eigenvalues $\lambda_1,\cdots,\lambda_m$ relies heavily on their population counterparts $\alpha_1,\cdots,\alpha_m$. 

Consider the functional inverse $\psi$ of the function $\alpha:~x\mapsto -1/\s(x)$. By \eqref{eq:st}, we have
\begin{equation}\label{eq:psi}
\psi(\alpha)=\psi_{y,H}(\alpha)=\alpha+y\alpha\int \frac{t}{\alpha-t}\d H(t),
\end{equation}
\[\psi'(\alpha)=1-y\int\frac{t^2}{(\alpha-t)^2}\d H(t),~\psi''(\alpha)=2y\int \frac{t^2}{(\alpha-t)^3}\d H(t).\]
This function $\psi(\cdot)$ is well defined for all $\alpha\notin \Gamma_H$. 

\cite{BY12} \hf{gave} a detailed characterization about the phase transition phenomenon of the limits of $\lambda_1,\cdots,\lambda_m$ when $\alpha_1,\cdots,\alpha_m$ satisfy different conditions. They name a generalized spiked eigenvalues $\alpha$ a {\em distant spike} for the M.P distribution $F^{y,H}$ if $\psi'(\alpha)>0$ and a {\em close spike} if $\psi'(\alpha)\leq 0$. Using the characterization of support of the LSD $F^{y,H}$ given in \cite{SilvChoi95}, it can be seen that for distant spikes, the corresponding sample eigenvalues almost surely converge to limits which are outside the support $\Gamma_{\underline{F}^{y,H}}$ of LSD of $\S_n$. These spikes are also referred as ``outliers" in the literature.

In this paper, we are focused on the generalized spiked model with distant spikes. \hf{In addition, for presentation simplicity, we only consider the case when $\alpha_1,\ldots,\alpha_m$ are non-identical. Extension to the case with possible overlaps on population spikes is straightforward using the developed techniques in this paper.} Therefore, in addition to Assumption\hf{s} (i) to (iv), we further assume that

\begin{itemize}
	\item[(v)] ${\bf \Lambda}$ is a fixed $m\times m$ matrix with non-identical bounded eigenvalues  $\alpha_1>\cdots>\alpha_m>\sup_p\max_j{\beta_{p',j}}$. All  $\alpha_k's$ are distant spiked eigenvalues satisfying $\psi'(\alpha_k)>0$.
\end{itemize}



\section{Main results}\label{sec:Joint}

In this section, we study the asymptotic behavior of $m$ largest sample  eigenvalues $\lambda_1,\cdots, \lambda_m$ and trace of $\S_n$.  Define the spectral decomposition of ${\bf \Lambda}=(\Lambda_{ij})_{m\times m}$ to be
\[
\mathbf{\Lambda}=\U ~\diag\lb \alpha_1,\cdots, \alpha_m
 \rb \U^{\T},\]
 where we remind that $\diag\lb\alpha_1,\cdots, \alpha_m\rb$ represents the $m\times m$ diagonal matrix with diagonal entries $\alpha_1,\cdots,\alpha_m$.

\begin{theorem}\label{thm:extreme}
	Under Assumptions (i) to (v), for the $m$ largest eigenvalues $\lambda_1,\cdots,\lambda_m$ of $\S_n$, denoting $\psi(\alpha_k)$ as $\psi_k$,  we have, the $m$-dimensional random vector
	$$\lb
	\sqrt{n}\lb \frac{\lambda_1}{\psi_1}-1\rb, \cdots,~
	\sqrt{n}\lb \frac{\lambda_m}{\psi_m}-1\rb
	\rb^\hf{\T}$$ and $\tr(\S_n)$ are jointly asymptotically normal and independent.
	Marginally,  
	\[\tr(\S_n)-\tr(\mathbf{\Sigma}_p)\xrightarrow{d}\EN(0, ~2y\gamma_2+y(\nu_4-3)\gamma_{d,2}),\]
	where for $k=1,2,\ldots$, $\gamma_k:=\int t^k \d H(t)$ denotes the $k$-th moment of LSD of $\V_{p'}$ and $\gamma_{d,2}:=\lim_{p'\rightarrow \infty}\frac{1}{p'}\sum_{i=1}^{p'}[\V_{p'}]_{ii}^2$.

	\noindent Moreover, $\lb
	\sqrt{n}\lb \frac{\lambda_1}{\psi_1}-1\rb, \cdots,~
	\sqrt{n}\lb \frac{\lambda_m}{\psi_m}-1\rb
	\rb^\hf{\T}$ weakly converges to an $m$-dimensional Gaussian vector 
	$\lb M_1,\cdots, M_m\rb^\hf{\T}$,  with each
	\[M_k=\bm{u}_k^{\T} \mathbf{G}(\psi(\alpha_k))\bm{u}_k,\]
	where $\bm{u}_k =\lb u_{1k},\cdots, u_{mk}\rb^{\T} $ is the $k$-th column of $\U $, $\mathbf{G}(\psi(\alpha_k))$ is \hf{an} $m\times m$ Gaussian random matrix with independent entries such that
	\begin{itemize}
		\item[(a)]its diagonal elements are i.i.d. Gaussian with mean zero and variance 
		\begin{equation}\label{eq:sigmak}
		\sigma^2_{\alpha_k}:=2\frac{\alpha_k^2}{\psi_k^2}\cdot\psi'(\alpha_k)+\beta_y\frac{\alpha_k^2}{\psi_k^2}\cdot(\psi'(\alpha_k))^2,
		\end{equation}
		with $\beta_y=\E|Y_{ij}|^4-3=\nu_4-3$ \hf{denoting} the fourth cumulant of base entries $\{Y_{ij}\}$; 
		\item[(b)] its upper triangular elements are i.i.d. Gaussian, with mean zero and variance
		\begin{equation}\label{eq:sk}
		s_{\alpha_k}^2=\frac{\alpha_k}{\psi_k}\cdot\psi'(\alpha_k).
		\end{equation}
	\end{itemize}
	Meanwhile, denoting the $(i,j)$-th entry of matrix $\A$ by $[\A]_{ij}$, we have, for $k_1\neq k_2$,
	{\small
\begin{gather*}
	\cov\lb \left[\G(\psi(\alpha_{k_1}))\right]_{ij}, \left[\G(\psi(\alpha_{k_2}))\right]_{ij}\rb=\frac{\alpha_{k_1}\alpha_{k_2}\psi'(\alpha_{k_1})\psi'(\alpha_{k_2})}{\psi_{k_1}\psi_{k_2}}\cdot\frac{\alpha_{k_1}-\alpha_{k_2}}{\psi_{k_1}-\psi_{k_2}},~i\neq j,\\
	\cov\lb \left[\G(\psi(\alpha_{k_1}))\right]_{ii}, \left[\G(\psi(\alpha_{k_2}))\right]_{ii}\rb=\frac{\alpha_{k_1}\alpha_{k_2}\psi'(\alpha_{k_1})\psi'(\alpha_{k_2})}{\psi_{k_1}\psi_{k_2}}\cdot\lb 2\cdot\frac{\alpha_{k_1}-\alpha_{k_2}}{\psi_{k_1}-\psi_{k_2}}+\beta_y\rb,\\
	\cov\lb \left[\G(\psi(\alpha_{k_1}))\right]_{ij}, \left[\G(\psi(\alpha_{k_2}))\right]_{i'j'}\rb=0, \mbox{ for all other cases.}
	\end{gather*}}
	\end{theorem}
Of note, \cite{WSY14} studied the limiting distribution of the random vector $$\lb \sqrt{n}\lb \frac{\lambda_1}{\psi_1}-1\rb, \cdots,~
	\sqrt{n}\lb \frac{\lambda_m}{\psi_m}-1\rb \rb^\hf{\T}$$ under Johnstone's spiked model \eqref{eq:Johnstone} \hf{with} $\V_{p'}=\I_{p'}$. Here we allow a more general $\V_{p'}$.

Following the proof of Theorem \ref{thm:extreme}, we \hf{are actually able to} provide a more accurate approximation for the asymptotic parameters of the limiting joint distribution of $(\lambda_k,\tr(\S_n))^{\T}$. Some second order terms which are ignored in the proof of Theorem \ref{thm:extreme} are sorted out in the following theorem. \hf{In applications, it may happen that the spiked eigenvalues $\alpha_k's$ are quite large while the sample size $n$ remains limited. In such situations, the correction terms below will be significant although their large $n$ limits are theoretically null. These terms are very useful for finite sample approximations as studied in the application of Section \ref{sec:App}. }

\begin{theorem}\label{prop:joint}
	 The asymptotic parameters in Theorem \ref{thm:extreme} can be further approximated as follows:
	{\small
	 \begin{gather*}
	\E \lb  \sqrt{n}\lb \frac{\lambda_k}{\psi_k}-1\rb\rb= \frac{\alpha_k^2\psi'(\alpha_k)}{\sqrt{n}\psi_k}\mu_{M_k}+ o\lb \frac{1}{\sqrt{n}\alpha_k^2}\rb,\\
	\var \lb  \sqrt{n}\lb \frac{\lambda_k}{\psi_k}-1\rb\rb =  \sum_{i=1}^mu_{ik}^4\sigma_{\alpha_k}^2+\sum_{i\neq j}^mu_{ik}^2u_{jk}^2s_{\alpha_k}^2+\frac{\alpha_k^4}{n\psi_k^2}(\psi'(\alpha_k))^2\sigma_{M_k}^2+o\lb\frac{1}{n\alpha_k^2}\rb,\\
	\var\lb\tr(\S_{n})-\tr(\mathbf{\Sigma}_{p})\rb= \frac{2}{n}\tr\lb \V_{p'}^2\rb+\frac{(\nu_4-3)}{n}\sum_{i=1}^{p'} [\V_{p'}]_{ii}^2+\frac{1}{n}\lb\sum_{i=1}^m\Lambda_{ii}^2(\nu_4-1)+\sum_{i\neq j} \Lambda_{ij}^2\rb,
	\end{gather*}}
\begin{equation*}
\resizebox{\textwidth}{!}{
	$\cov \lb  \sqrt{n}\lb \frac{\lambda_k}{\psi_k}-1\rb, \tr(\S_{n})-\tr(\mathbf{\Sigma}_{p})\rb= \rho_k+ \frac{y(\nu_4-1)}{\sqrt{n}\psi_k}\int\frac{t^2}{(1-t/\alpha_k)^2}\d H(t)+o\lb\frac{1}{\sqrt{n}\alpha_k}\rb$,}
\end{equation*}
	where 
	\begin{gather*}
	\rho_k=\frac{\alpha_k\psi'(\alpha_{k})}{\sqrt{n}\psi_k}\lb (\nu_4-1)\sum_{i=1}^m\Lambda_{ii}u_{ik}^2+\sum_{i\neq j}^m\Lambda_{ij}u_{ik}u_{jk}\rb, \\
	\mu_{M_k}=-\frac{1}{\alpha_k^3}\left(\cfrac{y\int\frac{t^2}{(1-t/\alpha_k)^3}\d H(t)}{\lb1-\frac{y}{\alpha_k^2}\int\frac{t^2}{(1-t/\alpha_k)^2}\d H(t)\rb^2}+\cfrac{y\beta_y\int\frac{t^2}{(1-t/\alpha_k)^3}\d H(t)}{1-\frac{y}{\alpha_k^2}\int\frac{t^2}{(1-t/\alpha_k)^2}\d H(t)}\right),\\
	\sigma_{M_k}^2=\cfrac{2\s'(\psi_k)\s'''(\psi_k)-3\lb\s''(\psi_k)\rb^2}{6\lb \s'(\psi_k)\rb^2}+y\beta_y\lb \s'(\psi_k)\rb^2\int\frac{t^2}{(1+t\s(\psi_k))^4}\d H(t).
	\end{gather*}
\end{theorem}

\begin{remark}
	\hf{In a related study, \cite{FW15} proved}, when the spiked part $\alpha_j=\alpha_j(p)\rightarrow \infty$ while  $c_j=p/(n\alpha_j)$ is bounded and the non-spiked part $$\frac{1}{p-m}\sum_{j=1}^{p'}\beta_{p',j}=\bar{c}+o(n^{-1/2}),$$ as $p/n\rightarrow \infty$, we have
	\[\sqrt{n}\left\{\frac{\lambda_j}{\alpha_j}-\lb1+\bar{c}c_j+O_{\mathbb{P}}(\alpha_j^{-1}\sqrt{p/n})\rb\right\}\xrightarrow{d} ~\EN(0,\nu_4-1).\]
	Notice that $\alpha_j^{-1}\sqrt{p/n}=c_j\sqrt{n/p}$, while $c_j$ is bounded and $p/n\rightarrow \infty$. \hf{Accordingly,} $O_{\mathbb{P}}(\alpha_j^{-1}\sqrt{p/n})$ is of order $o_{\mathbb{P}}(1)$. On the other hand, although our result in Theorem \ref{thm:extreme} is derived for bounded $\alpha_j$ for ease of presentation, it turns out that it is still valid when $\alpha_j\rightarrow \infty.$
	In fact, by some simple manipulations, Theorem \ref{thm:extreme} implies (under our settings) that
	\[\sqrt{n}\left\{\frac{\lambda_j}{\alpha_j}-\lb1+y\int \frac{t}{\alpha_j-t}\d H(t)\rb\right\}\xrightarrow{d} ~\EN(0,\nu_4-1).\]
It is \hf{interestingly compatible with} the result \hf{in} \cite{FW15} \hf{where} the term $\bar{c}c_j+O_{\mathbb{P}}(\alpha_j^{-1}\sqrt{p/n})$ is equivalent to $y\int \frac{t}{\alpha_j-t}\d H(t)$ \hf{as} $\alpha_j\rightarrow \infty$.	
\end{remark}

As an immediate application of Theorems \ref{thm:extreme} and \ref{prop:joint}, the following theorem characterizes the asymptotic behavior of an important statistic, the ratio statistic $\frac{\lambda_k}{p^{-1}\tr(\S_n)}$, which has been widely used in the  literature of signal detection.

\begin{theorem}\label{thm:ratio}
Under Assumptions (i) to (v), for $1\leq k\leq m$, we have
\begin{equation}
\sqrt{n}\lb\frac{\lambda_k}{\frac1p\tr(\S_n)}-\frac{\psi_k}{\frac1p\tr(\mathbf{\Sigma}_{p})}\rb\xrightarrow{d} \EN\lb 0,~\frac{\psi_k^2}{\gamma_1^2}\lb\sum_{i=1}^m u_{ik}^2\sigma_{\alpha_k}^2+\sum_{i\neq j}^m u_{ik}^2u_{jk}^2s_{\alpha_k}^2\rb\rb,
\end{equation}	
where $u_{ij}$ is the $(i,j)$-th entry of ${\bf \Lambda}$'s eigenmatrix $\U$,  $\gamma_1=\int t\d H(t)$, the mean value of LSD of $\V_{p'}$,  $\sigma_{\alpha_k}^2$ and $s_{\alpha_k}^2$ are defined in \eqref{eq:sigmak} and \eqref{eq:sk}.

\noindent
Moreover, the asymptotic variance in \eqref{thm:ratio} can be further approximated by
\begin{gather}\label{eq:ratiocorrect}
 \frac{\psi_k^2}{\lb\frac{1}{p}\tr(\mathbf{\Sigma}_p)\rb^2}\lb\sum_{i=1}^m u_{ik}^2\sigma_{\alpha_k}^2+\sum_{i\neq j}^m u_{ik}^2u_{jk}^2s_{\alpha_k}^2\rb 
+\frac{\alpha_k^4\lb\psi'(\alpha_k)\rb^2\sigma_{M_k}^2}{n\lb\frac{1}{p}\tr(\mathbf{\Sigma}_p)\rb^2}\\ \nonumber
-\frac{2\lb\psi_k+\frac{\alpha_k^2\psi'(\alpha_k)}{n}\mu_{M_k}\rb}{p\lb\frac{1}{p}\tr(\mathbf{\Sigma}_p)\rb^3}\left[\alpha_k\psi'(\alpha_k)\lb (\nu_4-1)\sum_{i=1}^m\Lambda_{ii}u_{ik}^2+\sum_{i\neq j}^m\Lambda_{ij}u_{ik}u_{jk}\rb\right.\\ \nonumber
\left.+y(\nu_4-1)\int\frac{t^2}{(1-t/\alpha_k)^2}\d H(t)\right]+\frac{\lb 2\gamma_2+(\nu_4-3)\gamma_{d,2}\rb\lb \psi_k+\frac{\alpha_k^2\psi'(\alpha_k)}{n}\mu_{M_k}\rb^2}{p\lb\frac{1}{p}\tr(\mathbf{\Sigma}_p)\rb^4}\\ \nonumber
+ \frac{\lb\sum_{i=1}^m\Lambda_{ii}^2(\nu_4-1)+\sum_{i\neq j} \Lambda_{ij}^2\rb\lb\psi_k+\frac{\alpha_k^2\psi'(\alpha_k)}{n}\mu_{M_k}\rb^2}{p^2\lb\frac{1}{p}\tr(\mathbf{\Sigma}_p)\rb^4}. \nonumber
\end{gather}
\end{theorem}
Under the generalized spiked population model in \eqref{eq:GeneralSpikeModel}, by applying Delta's method to the ratio function $f(x,y)=x/y$, with the joint limiting distribution of 
$\left(\lambda_1,\cdots,\lambda_m, \tr(\S_n)\right)^{\T}$ in Theorems \ref{thm:extreme} and \ref{prop:joint}, we immediately have the above limiting distribution of $\lambda_k/ \tr(\S_n)$ (proof thus omitted).
Observing that $p^{-1}\tr(\S_n)-p^{-1}\tr(\mathbf{\Sigma}_p)=o_p(1)$,  one might obtain \eqref{thm:ratio} by directly using Slutsky's theorem. However the second order terms in \eqref{eq:ratiocorrect} need implementation of the Delta's method.

\section{Application}\label{sec:App}

Determination of the number of signals embedded in noise is a fundamental problem in signal processing and chemometrics community. A significant portion of this literature has \hf{been} focused on the spiked covariance model \hf{a}rising from the following factor structure. Consider a sequence of $p$-dimensional random vectors $\{\x_t, ~t\in \mathbb{Z}\}$, admiting a version of static $m$-factor structure with $m$ fixed,
\begin{equation}\label{eq:factor}
\x_t=\A \bm{F}_t+\bm{E}_t.
\end{equation}
Here the factors $\bm{F}_t\sim \EN(0,\I_m)$ are assumed to be independent of the idiosyncratic error terms $\bm{E}_t\sim \EN(0,\sigma^2\I_p)$ with $\sigma^2$ fixed. The loading matrix $\A_{p\times m}$ is deterministic and of full rank such that $\A^{\T} \A $ has fixed eigenvalues $a_1>\cdots>a_m>0$. Suppose that we observe $\{\x_t,~t=1,\cdots,n\}$ with $n$ comparable to $p$.  In the high dimensional context, one major inference problem in \eqref{eq:factor} is to infer the total number of factors $m$.


Note that the eigenvalues of population covariance matrix $\mathbf{\Sigma}_{p}$ of $\x_t$ are
\[\mbox{Spec}\lb\mathbf{\Sigma}_p\rb=\{ a_1+\sigma^2,\cdots,a_m+\sigma^2,\underbrace{\sigma^2,\cdots,\sigma^2}_{p-m}\}.\]
Thus, \hf{it is immediate to observe that} $\x_t$ follows the generalized spiked model \eqref{eq:GeneralSpikeModel} with $\alpha_k=a_k+\sigma^2,~1\leq k\leq m$, $\beta_{p',1}=\cdots=\beta_{p',p'}=\sigma^2$. 

A typical \hf{approach to} test the number of factors is to find all the ``outliers" among eigenvalues $\lambda_1\geq\cdots\geq\lambda_p$ of sample covariance matrix.  According \hf{to} the phase transition phenomenon of sample eigenvalues established in \cite{BBP05}, \hf{it is known} that the asymptotic behavior of $\lambda_k~(1\leq k\leq m)$ differs depending on whether $\alpha_k/\sigma^2>1+\sqrt{p/n}$ or not.  Only when the signal-to-noise ratio (SNR) $\alpha_k/\sigma^2$ of the spikes is large enough can  the corresponding sample eigenvalues \hf{be separated away} from those bulk ones (outliers). Otherwise, these factors would be too weak and mix up with the noise.  In this section we develop a new test for presence of moderately strong factors. For any given integer $m_0\geq 1$ and constant $c> 1+\sqrt{p/n}$, we aim to test 
\begin{equation}\label{eq:factorH02}
H_0:~\frac{\alpha_{m_0}}{\sigma^2}\geq c,~\mbox{v.s. }~ H_1:~\frac{\alpha_{m_0}}{\sigma^2}< c.
\end{equation}
In other words, under $H_0$, there are at least $m_0$ signals with SNR larger than $c$. If $c=1+\sqrt{p/n}$, we are actually testing for the number of moderately strong factors above the critical transition value $1+\sqrt{p/n}$, that is to test
\begin{equation}\label{eq:factorH0}
H_0:~m\geq m_0, ~\mbox{v.s. }~ H_1:~ m< m_0.
\end{equation}
Compared to \cite{Choi17} and \cite{Onatski13} who have developed procedures to test ``$m=0$" against ``$m>0$" which still have non-zero power even when the factors are weak (SNR below the threshold $1+\sqrt{p/n}$), our test focuses more on testing existence of multiple strong spikes and their signal strength.

Since the noise variance $\sigma^2$ is typically unknown in real applications, we propose the following normalized sample eigenvalue as our test statistic,
$$T_{m_0}=\left.\lambda_{m_0}~\middle/~\frac{1}{p-m_0}\sum_{j>m_0}^p\lambda_j\right.$$ Here we use $(p-m_0)^{-1}\sum_{j>m_0}^p\lambda_j$ instead of the popular alternative $p^{-1}\tr(\S_n)$ as the surrogate for $\sigma^2$.  Although asymptotically equivalent under our conditions, 
$(p-m_0)^{-1}\sum_{j>m_0}^p\lambda_j$ is often found to have better finite-sample behavior and serve as a better estimate of $\sigma^2$ than $p^{-1}\tr(\S_n)$ under cases where several large spikes or a sizable collection of medium sized spikes are present. Correspondingly, $T_{m_0}$ often has superior detection power compared to $\widetilde{T}_{m_0}=\lambda_{m_0}/(p^{-1} \tr(\S_n))$ in finite sample cases. 

{ In the literature, \cite{Nadler11} and \cite{Deo16} adopted the test statistic $\widetilde{T}_1=\lambda_1/(p^{-1} \tr(\S_n))$ while they focus on testing the existence of one single spike, i.e.,
\begin{equation}\label{eq:Deotest}
H_0: \mathbf{\Sigma}_p=\sigma^2\I_p,~ \mbox{v.s. } ~H_1: m\geq 1.
\end{equation}
\noindent
Notice that this test is equivalent to the classical ``sphericity test with a spiked alternative". In the seminal paper \citep{Johnstone01}, it has been proven that when $\mathbf{\Sigma}_p=\sigma^2\I_p$  with Gaussian data, as $p,n\rightarrow \infty$ with $p/n\rightarrow y$,  
\[
\frac{\lambda_1/\sigma^2-\mu_{np}}{\sigma_{np}}\xrightarrow{d} {\rm TW}_{1},
\]
where
$\mu_{np}=\frac1n\lb\sqrt{n-1}+\sqrt{p}\rb^2,~\sigma_{np}=\sqrt{\frac{\mu_{np}}{n}}\lb \frac{1}{\sqrt{n-1}}+\frac{1}{\sqrt{p}}\rb^{1/3}$ and $\rm{TW}_{1}$ denotes Tracy-Widom distribution of order 1. \cite{Ma12} further refined the centering and scaling parameters
which improves the convergence rate from $O((n\wedge p)^{-1/3})$ to $O((n\wedge p)^{-2/3})$.
Note that the fluctuation in $\hat{\sigma}^2=p^{-1}\tr(\S_n)$
is of order $O(1/p)$, which is negligible compared to that in $\lambda_1$. Therefore, we still have,  as $p,n\rightarrow \infty$,
\[
\mathbb{P}\lb\frac{\widetilde{T}_1-\mu_{np}}{\sigma_{np}}<s\rb\rightarrow {\rm TW}_{1}(s).
\]
But in finite samples, the quality of this approximation for  $\widetilde{T}_1$ breaks down due to studentization. \cite{Nadler11} derived an explicit approximation formula for the tail probabilities of $\widetilde{T}_1$, which provides a significantly better fit to the empirical distribution of $\widetilde{T}_1$. But still, it is not a proper distribution function.
\cite{Deo16} suggested an alternative variance adjustment for the scaling parameter $\sigma_{np}$ to improve the finite sample performance of $\widetilde{T}_1$.
Tests based on these two corrections eliminate the downward size bias of the uncorrected test and improve its power performance for small values of $(p,n)$.


However, the problem of testing the existence of multiple spikes has not been fully resolved yet, i.e., for $m_0>1$, to test
\begin{equation}\label{eq:Deotestm}
H_0: m\leq m_0,~ \mbox{v.s. } ~H_1: m>m_0.
\end{equation}
First, in the presence of $m$ $(m>1)$ spikes, the conjecture of
\begin{equation}\label{eq:conjecture}
\cfrac{\lambda_{m+1}/\sigma^2-\mu_{n,p-m}}{\sigma_{n,p-m}}\xrightarrow{d}\rm{TW}_1
\end{equation}
has not been proven yet except for the complex Gaussian case (GUE case) in \cite{BBP05}. Secondly, the finite sample bias caused by the replacement of $\sigma^2$ with $p^{-1}\tr(\S_n)$ becomes more severe under $H_0$ {in \eqref{eq:Deotestm}}.  \cite{Deo16} and \cite{KN08} adopted alternative estimators of $\sigma^2$ based on $p-m$ bulk eigenvalues. However, their testing procedures are still based on the unverified conjecture \eqref{eq:conjecture}  and  simulation experiments show that these tests are still uniformly undersized due to the negative bias in the presence of spikes. Further analytic tools are needed to correct this bias.  }

In this paper we start from a different perspective by studying the behavior of {\red  $\lambda_m/\lb\frac{1}{p-m}\sum_{j>m}^p\lambda_j\rb$} instead of $\lambda_{m+1}$ in the presence of $m$ spikes and aim to test \eqref{eq:factorH02} and \eqref{eq:factorH0}. Although our test hypothesis is different from previous works, it can still be used to determine the total number of factors by performing a sequence of hypothesis tests on testing whether $\lambda_k$ $(1\leq k\leq m)$ arise from the signal or noise. 
The limiting distribution of our test statistic $T_{m}$ is fully implementable under $H_0$ in \eqref{eq:factorH02} assuming that all the conditions in Theorems \ref{thm:extreme} and \ref{prop:joint} are satisfied. Our test statistic is not only capable of testing the existence of multiple spikes, but \hf{can also} be used to test their signal strength. Higher order corrections are further made to alleviate finite sample bias, which ensures satisfactory testing size and power even when $(p,n)$ is not large. 

The corollary below {\red  follows from a direct implementation of Theorem \ref{prop:joint}}.

\begin{corollary}  
Considering the factor model \eqref{eq:factor}, for $1\leq k\leq m$, as $n\rightarrow\infty$, $p=p_n\rightarrow \infty$ such that $p/n\rightarrow y>0$, we have
\begin{equation}\label{eq:test}
\sqrt{n}\left(\frac{\lambda_k}{\frac{1}{p-k}\sum_{j>k}\lambda_j}-\frac{\psi_k}{\widetilde{\sigma}^2}\right)\xrightarrow{d} \EN\lb 0,~\frac{2\alpha_k^2\psi'(\alpha_k)}{\sigma^4}\rb,
\end{equation}
where \begin{gather*}
\psi_j=\alpha_j+\frac{y\alpha_j\sigma^2}{\alpha_j-\sigma^2},~\psi'(\alpha_j)=1-\frac{y\sigma^4}{(\alpha_j-\sigma^2)^2},~ \alpha_j=a_j+\sigma^2,\\
\widetilde{\sigma}^2=\frac{1}{p-k}\left[\tr\lb \Sigma_p\rb-\sum_{j=1}^k\left\{ \psi_j+\frac{\alpha_j^2\psi'(\alpha_j)}{n}\mu_{M_j} \right\}\right].
\end{gather*}
Additionally, we have a refined correction for the variance of $\frac{\lambda_k}{\frac{1}{p-k}\sum_{j>k}\lambda_j}$, 
\begin{align}\label{eq:factorcorr}
\sigma_{*,k}^2 =&\frac{2\alpha_k^2\psi'(\alpha_k)}{\widetilde{\sigma}^4}+\frac{4\alpha_k^2\psi'(\alpha_k)\left\{ \psi_k+\frac{2\alpha_k^2\psi_k'(\alpha_k)}{n}\mu_{M_k} \right\}}{(p-k)\widetilde{\sigma}^6}\\ \nonumber
&+\frac{\alpha_k^4\lb\psi'(\alpha_k)\rb^2\sigma_{M_k}^2}{n\widetilde{\sigma}^4}-\frac{4\alpha_k^2\left\{ \psi_k+\frac{\alpha_k^2\psi'(\alpha_k)}{n}\mu_{M_k} \right\}}{(p-k)\widetilde{\sigma}^6}\\ \nonumber
&+\frac{n}{(p-k)^2}\lb 2y\sigma^4+\frac{2}{n}\sum_{i=1}^m\alpha_i^2\rb\frac{\left\{ \psi_k+\frac{2\alpha_k^2\psi_k'(\alpha_k)}{n}\mu_{M_k} \right\}^2}{\widetilde{\sigma}^8}\\ \nonumber
&+\frac{2\psi_k\alpha_k^4 \lb\psi'(\alpha_k)\rb^2\sigma_{M_k}^2}{n(p-k)\widetilde{\sigma}^6}+\frac{\psi_k^2 \sum_{j=1}^k\left\{2\alpha_j^2\psi'(\alpha_j)-4\alpha_j^2\right\}}{(p-k)^2\widetilde{\sigma}^8}.
\end{align}
Here 
\begin{gather*}
\mu_{M_j}=-\frac{y\sigma^4}{\lb \alpha_j-\sigma^2\rb^3\left\{ 1-\frac{y\sigma^4}{\lb \alpha_j-\sigma^2\rb^2}\right\}^2}, ~\sigma_{M_j}^2=\cfrac{2\s'(\psi_j)\s'''(\psi_j)-3\lb\s''(\psi_j)\rb^2}{6\lb \s'(\psi_j)\rb^2}.
\end{gather*}
\end{corollary}


  Let $Z_{\alpha}$ be the lower-$\alpha$ quantile of the standard normal
  distribution at level $\alpha$.
  In order to define an appropriate critical value, one notes that
  the null hypothesis \eqref{eq:factorH02} is composite. For a given value of
  $\alpha_{m_0}$ under the null, the decision rule is to
  \[\mbox{reject } H_0 \text{ in} ~\eqref{eq:factorH02} {\text{ at the value } \alpha_{m_0}}~   \mbox{ if } T_{m_0}<\frac{\psi_{m_0}}{\widetilde{\sigma}^2}+Z_{\alpha}\cdot\frac{{\sigma_{*,m_0}}}{\sqrt{n}}.\]
  Letting  $t_k=\alpha_k/\sigma^2$, both the critical value above and the
  refined variance \eqref{eq:factorcorr} can be expressed as functions of
  these $t_k$'s:
  \begin{eqnarray*}
    \quad\quad q_{n,\alpha} &=& q_{n,\alpha}(t_{m_0},t_{m_0-1},\ldots,t_1)  =\frac{\psi_{m_0}}{\widetilde{\sigma}^2}+Z_{\alpha}\cdot\frac{{\sigma_{*,m_0}}}{\sqrt{n}}\nonumber \\
                 &=&\frac{ t_{m_0}+\frac{y}{1-1/t_{m_0}}}{1-\frac{1}{p-m_0}\sum_{j=1}^{m_0}\frac{y}{1-1/t_j}} +Z_{\alpha}\cdot\frac{\sigma_{*,m_0}}{\sqrt{n}},\nonumber
  \end{eqnarray*}
  \begin{align*}
{\rm and}~~~  \sigma_{*,m_0}^2=&  ~\sigma_{*,m_0}^2  (t_{m_0},t_{m_0-1},\ldots,t_1) \\ \nonumber
                    =& ~~2t_{m_0}^2\lb 1-\frac{y}{(t_{m_0}-1)^2}\rb\lb\frac{1}{1-\frac{1}{p-m_0}\sum_{j=1}^{m_0}\frac{y}{1-1/t_j}} \rb^2\\ \nonumber 
                    &-~\frac{4y t_{m_0}^2}{(p-m_0)(t_{m_0}-1)^2}\lb t_{m_0}+\frac{y}{1-1/t_{m_0}}\rb \lb\frac{1}{1-\frac{1}{p-m_0}\sum_{j=1}^{m_0}\frac{y}{1-1/t_j}} \rb^3\\ \nonumber   
                    &+~\frac{2yn}{(p-m_0)^2}\lb t_{m_0} +\frac{y}{1-1/t_{m_0}}\rb^2  \lb\frac{1}{1-\frac{1}{p-m_0}\sum_{j=1}^{m_0}\frac{y}{1-1/t_j}} \rb^4\\ \nonumber
                    &+~\frac{t_{m_0}^2\lb 1-\frac{y}{(t_{m_0}-1)^2}\rb^2}{n}\lb\frac{1}{1-\frac{1}{p-m_0}\sum_{j=1}^{m_0}\frac{y}{1-1/t_j}} \rb^2\left\{ \frac{4yt_{m_0}}{3(t_{m_0}-1)^3}\right.\\ \nonumber
                    &-~\frac{4yt_{m_0}}{3(t_{m_0}-1)^3\lb 1-\frac{y}{(t_{m_0}-1)^2}\rb^3} +\frac{2y^2t_{m_0}^2}{3(t_{m_0}-1)^6\lb 1-\frac{y}{(t_{m_0}-1)^2}\rb^4} \\ \nonumber
                    &+~\left.\frac{2y t_{m_0}^2}{(t_{m_0}-1)^4}+\frac{4y^2t_{m_0}^2}{3(t_{m_0}-1)^6\lb 1-\frac{y}{(t_{m_0}-1)^2}\rb}    \right\}   .
  \end{align*}
  Therefore, for the composite null \eqref{eq:factorH02}, we will use the
  critical value
  \begin{equation}
    \label{eq:qstar}
    q_{n,\alpha}^* =
    q_{n,\alpha}^*(t_{m_0-1},\ldots,t_1)=\inf_{c \le t_{m_0}< t_{m_0-1}}q_{n,\alpha}(t_{m_0},t_{m_0-1},\ldots,t_1),
  \end{equation}
  and reject the null if $T_{m_0}< q_{n,\alpha}^*$. For this procedure,
  it holds that
  \[  \limsup_{n\to\infty} \sup_{c \le t_{m_0}< t_{m_0-1}}\mathbb{P} (
    T_{m_0}< q_{n,\alpha}^*)\le \alpha.
  \] 

  Implementation of this procedure finally requires to estimate the
  values of larger spikes $\{t_k,~ 1\le k<m_0\}$ that appear in the
  critical value functions $q_{n,\alpha}$ and $q_{n,\alpha}^*$. As a matter of fact,
  consistent estimates for these spike values have been proposed in
  \cite{BD12}. 
  It is known that, under $H_0$ in \eqref{eq:factorH02}, for distant spikes $\alpha_k(1\leq k\leq  m_0)$, as $p,n\rightarrow \infty, ~p/n\rightarrow y$, $\lambda_k\xrightarrow{a.s.} \ell_k,$ satisfying $\s(\ell_k)=-\frac{1}{\alpha_k}$. Based on this result, \cite{BD12} proposed a consistent estimator $\hat{\alpha}_k$ for $\alpha_k$, 
  \[\hat{\alpha}_k=-\frac{1}{\hat{\s}(\ell_k)},\quad \hat{\s}(\ell_k)=-\frac{1-y}{\lambda_k}+\frac{1}{n}\sum_{i>m_0}\frac{1}{\lambda_i-\lambda_k}.\]
  Here  the noise level $\sigma^2$ is intrinsically hidden inside the values of sample eigenvalues $\lambda_j$'s, not explicitly used in the estimation of $\alpha_k$. Note that the conditions used in \cite{BD12} which ensure the consistency of $\hat{\alpha}_k$ are satisfied under our settings (Assumptions (i) to (v)).
  Thus for $1\leq k<m_0$, the  $t_k$'s can be consistently estimated  by
  \[\hat{t}_k=\frac{\hat{\alpha}_k}{\lb\tr(\S_n)-\sum_{\ell=1}^{m_0}\hat{\alpha}_{\ell}\rb/(p-m_0)}.\]
  Plugging these estimates into the critical value functions
  $q_{n,\alpha}$ and $q_{n,\alpha}^*$ leads to  a
  full implementation of the test.

  To summarize, the proposed testing procedure is to 
  \begin{equation}\label{eq:rejectreal}
    \mbox{reject } H_0\text{ in }\eqref{eq:factorH02} \mbox{ if } T_{m_0}< {\hat q}_{n,\alpha}^*, 
  \end{equation}
  where
  \begin{equation}
    \label{eq:hatq}
    {\hat q}_{n,\alpha}^* =q_{n,\alpha}^*(\hat t_{m_0-1},\ldots,\hat t_1)
    = \inf_{c \le t_{m_0}< \hat
      t_{m_0-1}}q_{n,\alpha}(t_{m_0},\hat t_{m_0-1},\ldots,\hat t_1).
  \end{equation}

{Simulation experiments are conducted in the supplement to examine the performance of out testing procedure.  Empirical data are generated following the factor model in
\eqref{eq:factor}.  Different numbers of factors and signal strength are assigned for various model settings. Results (Tables \ref{tab:size} and \ref{tab:power}) show that our test works for different numbers of factors.  It has  better
performance for higher signal strength levels. The same supplement section has also designed experiments to illustrate the necessity of incorporating the second order correction to the asymptotic variance  proposed  in \eqref{eq:factorcorr}. Numerical  comparison is made
  between \eqref{eq:rejectreal} and the following testing procedure to
  \begin{align}\label{eq:reject2}
    &\mbox{reject } H_0\text{ in }\eqref{eq:factorH02} \mbox{ if } T_{m_0}<  {\tilde q}_{n,\alpha}^*,
  \end{align}
where ${\tilde q}_{n,\alpha}^*$ is defined as ${\hat
    q}_{n,\alpha}^*$ in \eqref{eq:hatq} except that $\frac{\psi_{m_0}}{\widetilde{\sigma}^2}$   is replaced by $\frac{\psi_{m_0}}{\sigma^2}$ and the refined asymptotic
  variance $\sigma_{*,m_0}^2$ used there is replaced by
  $
    \tilde{\sigma}^2_{m_0}= 2t_{m_0}^2-\frac{2yt_{m_0}^2}{(t_{m_0}-1)^2}.
  $
Tables \ref{tab:corr} and \ref{tab:corr2} in the supplement show that  this refined
correction \eqref{eq:factorcorr} for the  variance plays an important role in controlling the size and improves the power in the testing procedure, especially when the data dimension and sample size are relatively small. }

\section{Proofs}\label{sec:proofs}

\subsection{Proof of Theorem \ref{thm:extreme}}

Considering the block structure of population covariance matrix $\mathbf{\Sigma}_p$, the analysis is carried out using a decomposition into blocks of size $m$ and $p-m$ respectively. Define 
\[
\y_i=\lb
\begin{array}{c}
\y_{1i}\\
\y_{2i}
\end{array}
\rb, \quad\mathbf{\Sigma}_p=\lb
\begin{array}{cc}
{\bf\Lambda}&{\bf 0}\\
{\bf 0} & {\bf V}_{p'}
\end{array}\rb,\quad\x_i=\lb 
\begin{array}{c}
\x_{1i}\\
\x_{2i}
\end{array}
\rb=\lb\begin{array}{c}
{\bf \Lambda}^{1/2}\y_{1i}\\
\V_{p'}^{1/2}\y_{2i}
\end{array}
\rb.~
\]
The sample covariance matrix is then
\[\S_n=\frac{1}{n}\sum_{i=1}^n\x_i\x_i^{\T}  =\frac{1}{n}\lb
\begin{array}{cc}
\X_1\X_1^{\T}   & \X_1\X_2^{\T}  \\
\X_2\X_1^{\T}   & \X_2\X_2^{\T}  
\end{array}
\rb=\lb
\begin{array}{cc}
\S_{11} & \S_{12}\\
\S_{21} & \S_{22}
\end{array}
\rb,
\]
where 
$$\S_{11}=\frac{1}{n}{\bf \Lambda}^{1/2}\Y_1\Y_1^{\T}  {\bf \Lambda}^{1/2},\quad\S_{22}=\frac{1}{n}{\V_{p'}}^{1/2}\Y_2\Y_2^{\T}  {\V_{p'}}^{1/2},$$
$$\X_1=\lb \x_{11},\cdots,\x_{1n}\rb,~ \X_2=\lb\x_{21},\cdots,\x_{2n}\rb,$$
$$\Y_1=\lb \y_{11},\cdots,\y_{1n}\rb,~ \Y_2=\lb\y_{21},\cdots,\y_{2n}\rb.$$

The proof of Theorem \ref{thm:extreme} builds on the above block-decomposition analysis of spiked and non-spiked ones. $\tr(\S_n)$ follows the decomposition $\tr(\S_n)=\tr(\S_{11})+\tr(\S_{22})$. It will later be shown that $\tr(\S_{22})$ is asymptotically independent of the random vector $\lb\lambda_1,\cdots,\lambda_m\rb^{\T}$\hf{, while the covariance between $\tr(\S_{11})$ and $\lb\lambda_1,\cdots,\lambda_m\rb^{\T}$ is of the order $O(1/\sqrt{n})$.} 
The proof in general consists of four steps as follows:

\begin{description}
	\item[Step 1.] \hf{deriving the} asymptotic joint distribution of $(\lambda_1,\cdots,\lambda_m)^{\T}$;
	\item[Step 2.] \hf{deriving the} marginal limiting distribution of $\tr(\S_{11})$ and $\tr(\S_{22})$;
	\item[Step 3.] \hf{deriving the} asymptotic joint distribution of $\lb \lb\lambda_k\rb_{1\leq k\leq m},\tr(
	\S_{n})\rb^{\T}$.
\end{description}

\medskip\noindent{\bf \textsc{Step 1: Joint limiting distribution of $\lb\lambda_1,\cdots,\lambda_m\rb^{\T}$.}}

\noindent
Many efforts in the literature have been put into the study of the asymptotic behavior of extreme sample eigenvalues under various spiked population models. Notably, \cite{BY12} derived a CLT for sample eigenvalues correspond\hf{ing} to one distant spiked eigenvalue under a more generalized model where possible multiplicity of spiked eigenvalues is allowed such that
\begin{equation}\label{eq:GeneralSpikeModel2}
\mbox{Spec}(\mathbf{\Sigma}_p)=\Big\{\underbrace{\alpha_1,\cdots,\alpha_1}_{m_1},\cdots,\underbrace{\alpha_K,\cdots,\alpha_K}_{m_K},~\beta_{p',1},\cdots,\beta_{p',p'}\Big\}.
\end{equation}

Here we eliminate the multiplicity of spikes in model \eqref{eq:GeneralSpikeModel2} and focus on the correlation among sample eigenvalues corresponding to different spikes. \cite{WSY14} studied such correlations under the model where 
\[\mbox{Spec}(\mathbf{\Sigma}_p)=\Big\{\underbrace{\alpha_1,\cdots,\alpha_1}_{m_1},\cdots,\underbrace{\alpha_K,\cdots,\alpha_K}_{m_K},1,\cdots,1\Big\},\]
\hf{while} our model \eqref{eq:GeneralSpikeModel} is more general in the sense of bulk eigenvalue distribution $H(t)$. The proof in general combines the Z-estimation scheme in \cite{BY12} and the result of joint CLT for several random sesquilinear forms in \cite{WSY14}. Detailed proofs are presented here, which will also be used in subsequent steps and the proof of Theorem \ref{prop:joint}.

Noting that,  whenever $\A$ is invertible,
\[
\left\lvert 
\begin{array}{cc}
\A & \mathbf{B}\\
\mathbf{C} & \mathbf{D}
\end{array}
\right\rvert=|\A|\cdot|\D-\C\A^{-1}\B|,
\]
an eigenvalue $\lambda_i$ of $\S_n$ that is not an eigenvalue of $\S_{22}$ satisfies
\begin{equation}\label{eq:decompose}
0=\left\lvert \lambda_i \I_p-\S_n\right\rvert=\left\lvert \lambda_i \I_{p-m}-\S_{22}\right\rvert\cdot\left\lvert \lambda_i \I_m-\K_n(\lambda_i)\right\rvert,
\end{equation}
where 
\[\K_n(\ell):=\S_{11}+\S_{12}\lb \ell\I_{p-m}-\S_{22}\rb^{-1}\S_{21}.
\]
Thus, the eigenvalues of $\S_n$ satisfy
\begin{equation}\label{eq:lambdai}
\left\lvert \lambda_i \I_m-\K_n(\lambda_i)\right\rvert=0.
\end{equation}

Consider a real number $\ell$ outside the support of LSD $F^{y,H}$ of $\S_{22}$ and the goal is to find the limit of random matrix $\K_n(\ell)$ with fixed dimension $m$. Since for $\ell\neq 0$ which is not an eigenvalue of $\A^{\T}  \A$,
\[\I_n+\A\lb \ell\I_{p-m}-\A^{\T}  \A\rb^{-1}\A^{\T}  =\ell\lb \ell\I_n-\A\A^{\T}  \rb^{-1},\]
it holds that
\begin{align}
\K_n(\ell)~&=~\frac{1}{n}\X_1\left\{\I_n+\frac{1}{n}\X_2^{\T}  \lb \ell\I_{p-m}-\frac{1}{n}\X_2\X_2^{\T}   \rb^{-1}\X_2\right\}\X_1^{\T}  \\ \nonumber
~&=~\frac{\ell}{n} {\bf \Lambda}^{1/2}\Y_1\lb \ell\I_n-\frac{1}{n}\Y_2^{\T}  \V_{p'}\Y_2\rb^{-1}\Y_1^{\T}  {\bf \Lambda}^{1/2}.
\end{align}

By Assumptions (i)-(iv), $\ell$ is outside the support of LSD $F^{y,H}$ of $\S_{22}$ and for $n$ large enough, the operator norm of $\lb \ell\I_{p-m}-\frac{1}{n}\Y_2^{\T}  \V_{p'}\Y_2 \rb^{-1}$ is bounded. Meanwhile, $\Y_1$ is independent of $\Y_2$. Then by the law of large numbers, Corollary 3.9 in \cite{KY14},  and Theorem 11.8 in \cite{YZB15a}, we have, as $p,n\rightarrow \infty$, $p/n\rightarrow y>0$,

\[{\bf\Lambda}^{-1/2} \K_n(\ell){\bf\Lambda}^{-1/2}\xrightarrow{a.s.}-\ell\s(\ell)\cdot\I_m,\]
where $\s(\ell)$ is the Stieltjes transform of the LSD of $\frac{1}{n}\Y_2^{\T}  \V_{p'}\Y_2$.

Therefore, the eigenvalue $\lambda_i$ of $\S_{n}$ satisfying \eqref{eq:lambdai} converges almost surely to \hf{the limit} $\ell$ such that
\[\left\lvert \ell\cdot{\bf \Lambda}^{-1}+\ell\s(\ell)\cdot\I_m\right\rvert=0,\]
that is 
\[\s(\ell)=-1/\alpha_k,~\ell=\psi(\alpha_k),~ k=1,\cdots,m,\]
where $\mbox{Spec}({\bf \Lambda})=\{\alpha_1,\cdots,\alpha_m\}$. 
The following lemma, due to \cite{SilvChoi95}, characterizes the close relationship between the supports of the generating measure $H$ and the \hf{Mar\v{c}enko-Pastur (M-P)} distribution $F^{y,H}$.

\begin{lemma}\label{prop:support}(By \cite{SilvChoi95})
	If $\lambda\notin \Gamma_{\underline{F}^{y,H}}$, then $\s(\lambda)\neq 0$ and $\alpha=-1/\s(\lambda)$ satisfies
	\begin{itemize}
		\item[(1)] $\alpha\notin \Gamma_H$ and $\alpha\neq 0$ (so that $\psi(\alpha)$ is well defined);
		\item[(2)] $\psi'(\alpha)>0$.
	\end{itemize}
	Conversely, if $\alpha$ satisfies (1)-(2), then $\lambda=\psi(\alpha)\notin \Gamma_{\underline{F}^{y,H}}$.
\end{lemma}

By Lemma \ref{prop:support}, $\ell=\psi(\alpha_k)$ is outside the support of LSD $F^{y,H}$
if and only if $\psi'(\alpha_k)>0$. By Assumption (v), all spiked values $\alpha_k$ are large enough to make $\psi'(\alpha_k)>0$. Therefore, the limits 
$$\ell=\psi(\alpha_k)=:\psi_k, ~k=1,\cdots,m$$ 
are all outside the support of  LSD $F^{y,H}$. Meanwhile, for the $m$ largest eigenvalues $\lambda_1,\cdots,\lambda_m$ of $\S_n$, as $p,n\rightarrow\infty$, $p/n\rightarrow y>0$,
\[\frac{\lambda_k}{\psi_k}\xrightarrow{a.s.} 1,~ \mbox{ for } 1\leq k\leq m.\]
Note that $\psi(\cdot)$ is the functional inverse of function $\alpha:~x\mapsto -1/\s(x)$, then

\begin{equation}\label{eq:s}
\s(\psi_k)=-\frac{1}{\alpha_k}, ~ \s'(\psi_k)=\frac{1}{\alpha_k^2\psi'(\alpha_k)},~\s''(\psi_k)=-\frac{2}{\alpha_k^3(\psi'(\alpha_k))^2}-\frac{\psi''(\alpha_k)}{\alpha_k^2(\psi'(\alpha_k))^3},
\end{equation}

\begin{align*}
\s'''(\psi_k)&=\frac{6}{\alpha_k^4(\psi'(\alpha_k))^3}-\frac{2\psi''(\alpha_k)\s''(\psi_k)}{(\psi'(\alpha_k))^2}\\
&-\frac{\alpha_k^2\psi'(\alpha_k)\psi'''(\alpha_k)-\psi''(\alpha_k)\lb 2\alpha_k\psi'(\alpha_k)+\alpha_k^2\psi''(\alpha_k)\rb}{\alpha_k^4(\psi'(\alpha_k))^2}, \nonumber
\end{align*}
\noindent
where $\left.[\cdot]'\right\rvert_{\alpha=\alpha_k}$ means to take first order derivative with regard to $\alpha$ and then let $\alpha=\alpha_k$. Here
\begin{gather*}
\psi_k=\psi(\alpha_k)=\alpha_k+y\int \frac{t}{1-t/\alpha_k}\d H(t),\\
\psi'(\alpha_k)=1-y\int\frac{t^2}{(\alpha_k-t)^2}\d H(t),~\psi''(\alpha_k)=2y\int\frac{t^2}{(\alpha_k-t)^3}\d H(t).
\end{gather*}

Denote  $\s_n(\ell)=\frac{1}{n}\tr\lb\ell\I_n-\frac{1}{n}\Y_2^{\T}\V_{p'}\Y_2\rb^{-1}$, then
\begin{align*}
\K_{n}(\ell) 
~& = ~ 
 \frac{\ell}{\sqrt{n}}\mathbf{\Lambda}^{1/2}\R_{n}\left(\ell\right)\mathbf{\Lambda}^{1/2}-\ell\mathbf{\Lambda}\s(\ell)+\ell\mathbf{\Lambda}\left(\s(\ell)-\s_{n}(\ell)\right),
\end{align*}
where $\R_n(\ell)$ is a sequence of $m\times m$ random matrix-valued processes
\begin{equation}\label{eq:Rnl}
\resizebox{\textwidth}{!}{
$\R_n(\ell)=\frac{1}{\sqrt{n}}\lb \Y_1\lb \ell\I_n-\frac{1}{n}\Y_2^{\T}  \V_{p'}\Y_2\rb^{-1}\Y_1^{\T}  -\I_m\tr\lb \ell\I_n-\frac{1}{n}\Y_2^{\T}  \V_{p'}\Y_2\rb^{-1}\rb,~\ell\in \mathcal{U}.$}
\end{equation}
\hf{Here} $\mathcal{U}$ is a compact set of indexes outside the support of LSD of $\S_{22}$. 

The establishment of CLT for extreme sample eigenvalues $\lambda_i's$ relies heavily on the finite dimensional convergence of processes \[\left\{\R_n(\ell),~\ell\in \mathcal{U}\right\}~\mbox{and }~\left\{n(\s_n(z)-\s(z)),~z\in \mathbb{C}\setminus \Gamma_{F^{y,H}}\right\},\]
which has been well established in \cite{BY12}, Lemma 1.1 in \cite{BS04}. More specifically, we have the following lemma.

\begin{lemma}\label{prop:Rn}(By Theorem 11.10 in \cite{YZB15a})
Under Assumptions (i) to (iv), for any $L$ index values $\{\ell_j\}$, the  $L$ random matrices 
\[\Big\{ \R_n(\ell_1),\cdots,\R_n(\ell_L)\Big\}\]
 weakly converge to $L$ Gaussian random matrices determined as follows: for arbitrary $L$ numbers $a_1,\cdots, a_L$, the random matrix
\[\widetilde{\R}_n=a_1\R_n(\ell_1)+\cdots+a_L\R_n(\ell_L)\]
weakly converges to a Gaussian random matrix $\R=\{R_{ij}\}$ where
\begin{itemize}
	\item[(1)] the diagonal entries are i.i.d. zero-mean Gaussian with variance 
	\[\var(R_{ii})=w(\E|Y_{ij}|^4-3)+2\theta;\]
	\item[(2)] the upper off-diagonal entries are i.i.d. zero-mean Gaussian with variance $\theta$;
	\item[(3)] all these entries are mutually independent.
\end{itemize}
Here the parameters $\theta$ and $w$ are 
\[\theta=\sum_{j=1}^La_j^2\s'(\ell_j)+2\sum_{j<k}a_ja_k\frac{\s(\ell_j)-\s(\ell_k)}{\ell_j-\ell_k},~w=\lb\sum_{j=1}^La_j\s(\ell_j)\rb^2.\]
\end{lemma}

\noindent
Notice that
\[
\lambda_{k}\I_{m}-\K_{n}(\lambda_{k})=\lambda_{k}\I_{m}-\frac{\lambda_{k}}{\sqrt{n}}\mathbf{\Lambda}^{1/2}\R_{n}(\lambda_{k})\mathbf{\Lambda}^{1/2}+\lambda_{k}\mathbf{\Lambda}\s(\lambda_{k})-\lambda_{k}\mathbf{\Lambda}\left(\s(\lambda_{k})-\s_{n}(\lambda_{k})\right).
\]

\noindent
Since we have spectral decomposition
${\bf \Lambda}=\U \diag\lb \alpha_1 ,\cdots, \alpha_m\rb \U^{\T}  $,
\begin{eqnarray*}
&~&\U^{\T}\left(\I_{m}-\lambda_{k}^{-1}\K_{n}(\lambda_{k})\right)\U =\left(\begin{array}{ccc}
	\ddots\\
	& (1+\alpha_{u}\s(\lambda_{k}))\\
	&  & \ddots
\end{array}\right)_{m\times m}\\
&-&\frac{1}{\sqrt{n}}\U^{\T}\mathbf{\Lambda}^{1/2}\R_{n}(\lambda_{k})\mathbf{\Lambda}^{1/2}\U-\left(\begin{array}{ccc}
\ddots\\
& \alpha_{u}\left(\s(\lambda_{k})-\s_{n}(\lambda_{k})\right)\\
&  & \ddots
\end{array}\right)_{m\times m}.
\end{eqnarray*}

\noindent
Now consider\hf{ing} $\delta_{n,k}=\sqrt{n}\left(\frac{\lambda_{k}}{\psi_{k}}-1\right)$, $\psi_{k}=\alpha_{k}+y\int\frac{t}{1-t/\alpha_{k}}\d H(t)$,
by Taylor expansion, we have
\[
\s(\lambda_{k})=\s(\psi_{k})+\frac{\psi_{k}}{\sqrt{n}}\cdot\delta_{n,k}\cdot\s'(\psi_k)+O_{\mathbb{P}}\lb\frac{1}{n}\rb,
\]
which then yields
\begin{eqnarray*}
	&~&\U^{\T}\left(\I_{m}-\lambda_{k}^{-1}\K_{n}(\lambda_{k})\right)\U\\
	 & = & \left(\begin{array}{ccc}
		\ddots\\
		& 1+\alpha_{u}\s(\psi_{k})+\alpha_{u}\lb\frac{\psi_{k}}{\sqrt{n}}\cdot\delta_{n,k}\cdot\s'(\psi_k)+O_{\mathbb{P}}\lb\frac{1}{n}\rb\rb\\
		&  & \ddots
	\end{array}\right)\\
	&& -  \frac{1}{\sqrt{n}}\U^{\T}\mathbf{\Lambda}^{1/2}\R_{n}(\lambda_{k})\mathbf{\Lambda}^{1/2}\U-\left(\begin{array}{ccc}
		\ddots\\
		& \alpha_{u}\left(\s(\lambda_{k})-\s_{n}(\lambda_{k})\right)\\
		&  & \ddots
	\end{array}\right).
\end{eqnarray*}

First, it can be seen that all the \hf{non-}diagonal terms tend to zero
on the right hand side. Then for a diagonal term with index $u\neq k$,
by definition $1+\s(\psi_{k})\alpha_{u}\neq0$ and it is indeed the
leading term with the remaining three terms converg\hf{ing}
to zero. As for the $k$-th diagonal term, $1+\s(\psi_{k})\alpha_{k}=0$
by definition, thus the $k$-th diagonal term reduces to 
\[
\resizebox{\textwidth}{!}{$
\alpha_{k}\lb\frac{\psi_{k}}{\sqrt{n}}\cdot\delta_{n,k}\cdot\s'(\psi_k)+O_{\mathbb{P}}\lb\frac{1}{n}\rb\rb-\frac{1}{\sqrt{n}}\left[\U^{\T}\mathbf{\Lambda}^{1/2}\R_{n}(\lambda_{k})\mathbf{\Lambda}^{1/2}\U\right]_{kk}-\alpha_{k}\left(\s(\lambda_{k})-\s_{n}(\lambda_{k})\right).$}
\]

\noindent
Noting that 
$\left\lvert\U^{\T}\left(\I_{m}-\lambda_{k}^{-1}\K_{n}(\lambda_{k})\right)\U\right\rvert=0$,  for $n$ sufficiently large,  we have
\begin{equation}\label{eq:delta}
\left\vert \frac{\psi_{k}}{\sqrt{n}}\cdot\delta_{n,k}\cdot\s'(\psi_k)+O_{\mathbb{P}}\lb\frac{1}{n}\rb-\frac{1}{\sqrt{n}}\left[\U^{\T}\R_{n}(\lambda_{k})\U\right]_{kk}-\left(\s(\lambda_{k})-\s_{n}(\lambda_{k})\right)\right\rvert=0.
\end{equation}
Taking into account the convergence of process $\left\{ \R_{n}(\ell),~\ell\in\mathcal{U}\right\} $  and \\$\left\{ M_{n}(z)=n\lb\s_{n}(z)-\s(z)\rb\right\}$ in Lemma \ref{prop:Rn} and \cite{BS04}, 
it follows that $\delta_{n,k}$ weakly converges to a solution of the limit of

\[
\left\lvert \frac{\psi_{k}}{\sqrt{n}}\cdot\delta_{n,k}\cdot\s'(\psi_{k})-\frac{1}{\sqrt{n}}\left[\U^{\T}\R_{n}(\psi_{k})\U\right]_{kk}-O_{\mathbb{P}}\lb\frac{1}{n}\rb\right\rvert =0,
\]
i.e.,
\[
\delta_{n,k}=\sqrt{n}\left(\frac{\lambda_{k}}{\psi_{k}}-1\right)\xrightarrow{d} \frac{\left[\U^{\T}\R(\psi_{k})\U\right]_{kk}}{\psi_{k}\s'(\psi_{k})}.
\]

\noindent
Denote $$\G(\psi_k)=\frac{\R(\psi_k)}{\psi_k\s'(\psi_k)}.$$ Then by Lemma \ref{prop:Rn}, 
$\G(\psi_k)$ is a Gaussian random matrix with mutually independent entries where the diagonal entries are i.i.d. zero-mean Gaussian with variance
\[\var(\left[\G(\psi_k)\right]_{ii})=\cfrac{(\nu_4-3)\s(\psi_k)^2+2\s'(\psi_k)}{(\s'(\psi_k))^2\psi_k^2},\]
and the upper off-diagonal entries are i.i.d. zero mean Gaussian with variance 
\[\var(\left[\G(\psi_k)\right]_{ij})=\frac{1}{\s'(\psi_k)\psi_k^2}.\]
\noindent

Now we consider the asymptotic joint distribution of $(\lambda_{k_1}, \lambda_{k_2})^{\T},~1\leq k_1\neq k_2\leq m$. It can be seen from \hf{the} previous proof that the leading term of $\delta_{n,k}=\sqrt{n}\lb\frac{\lambda_k}{\psi_k}-1\rb$ is $\frac{\left[\U^{\T}  \R_n(\psi_k)\U\right]_{kk}}{\s'(\psi_k)\psi_k}$. Thus the correlation between limits of $(\lambda_{k_1},\lambda_{k_2})^{\T}$ is determined by the joint limiting distribution of the two random sesquilinear forms $\R_n(\psi_{k_1})$ and $\R_n(\psi_{k_2})$. \hf{This task is nontrivial.} 
Here we apply a joint CLT for random vector whose components are function of random sesquilinear forms by \cite{WZY14} .

\begin{lemma}\label{lem:jointRn} (By \cite{WZY14})
	Consider a sequence $(\x_i,\y_i)_{i\in \mathbb{N}}$ of i.i.d. real valued zero-mean random vector\hf{s} belonging to $\mathbb{R}^K\times\mathbb{R}^K$ with finite fourth order moment\hf{:}
	\begin{gather*}
	\x_i=(X_{1i},\cdots,X_{Ki})^{\T},1\leq i\leq n,~  \x(\ell )=(X_{\ell 1},\cdots,X_{\ell n})^{\T},1\leq \ell\leq K,\\
	\y_i=(Y_{1i},\cdots,Y_{Ki})^{\T},1\leq i\leq n,~  \y(\ell)=(Y_{\ell 1},\cdots,Y_{\ell n})^{\T},1\leq \ell\leq K,
	\end{gather*}
	\hf{and} $\rho(\ell)\hf{:}=\mathbb{E}(X_{\ell 1}Y_{\ell 1})$. Let $\{\A_n\}$ and $\{\B_n\}$ be two sequence\hf{s} of $n\times n$ symmetric matrices. Assume the following limits exist:
	\begin{gather*}
	w_1=\lim_{n\rightarrow \infty}\frac1n\tr(\A_n\odot \A_n),~ w_2=\lim_{n\rightarrow \infty}\frac1n\tr(\B_n\odot \B_n),~w_3=\lim_{n\rightarrow \infty}\frac1n\tr(\A_n\odot \B_n),\\
	\theta_1=\lim_{n\rightarrow \infty}\frac1n\tr(\A_n \A_n^{\T}  ),~
	\theta_2=\lim_{n\rightarrow \infty}\frac1n\tr(\B_n\B_n^{\T}  ),~
	\theta_3=\lim_{n\rightarrow \infty}\frac1n\tr(\A_n\B_n^{\T}  ),
	\end{gather*}
	where $\A\odot \B$ denotes the Hadamard product of two matrices $\A$ and $\B$, i.e., $[\A\odot \B]_{ij}\hf{:}=[\A]_{ij}[\B]_{ij}$. Define two groups of sesquilinear forms:
{\small
	\[U(\ell)=\frac{1}{\sqrt{n}}\lb \x(\ell)^{\T}  \A_n\y(\ell)-\rho(\ell)\tr(\A_n)\rb,~V(\ell)=\frac{1}{\sqrt{n}}\lb \x(\ell)^{\T}  \B_n\y(\ell)-\rho(\ell)\tr(\B_n)\rb.\]}
	Then the $2K$-dimensional random vector $(U(1),\cdots,U(K),V(1),\cdots,V(K))^\hf{\T}$ weakly converges to a zero-mean Gaussian vector with covariance matrix $\B=\lb
	\begin{array}{cc}
	\B_{11}&\B_{12}\\
	\B_{21}&\B_{22}
	\end{array}
	\rb_{2K\times 2K}$, with each block $\B_{ij}=\lb\B_{ij}(\ell,\ell')\rb_{1\leq \ell,\ell'\leq K}$ a $K\times K$ matrix having structure, for $1\leq \ell,\ell'\leq K$,
	\begin{gather*}
	\B_{11}(\ell,\ell')=\cov(U(\ell),U(\ell'))=w_1a_1+(\theta_1-w_1)(a_2+a_3),\\
	\B_{22}(\ell,\ell')=\cov(V(\ell),V(\ell'))=w_2a_1+(\theta_2-w_2)(a_2+a_3),\\
\B_{12}(\ell,\ell')=\cov(U(\ell),V(\ell'))=w_3a_1+(\theta_3-w_3)(a_2+a_3),\\
	a_1=\mathbb{E}(X_{\ell 1}Y_{\ell 1}X_{\ell'1}Y_{\ell'1})-\rho(\ell)\rho(\ell'),\\
	~a_2=\mathbb{E}(X_{\ell 1}X_{\ell'1})\mathbb{E}(Y_{\ell 1}Y_{\ell'1}),~a_3=\mathbb{E}(X_{\ell 1}Y_{\ell'1})\mathbb{E}(X_{\ell'1}Y_{\ell 1}).
	\end{gather*}
\end{lemma}

\noindent
Not\hf{ing} that 
{\small
\[\R_n(\ell)=\frac{1}{\sqrt{n}}\lb \Y_1\lb \ell\I_n-\frac{1}{n}\Y_2^{\T}  \V_{p'}\Y_2\rb^{-1}\Y_1^{\T}  -\I_m\tr\lb \ell\I_n-\frac{1}{n}\Y_2^{\T}  \V_{p'}\Y_2\rb^{-1}\rb,
\]}
it can be seen that, for any $1\leq i,j,i',j' \leq m$, the  random vector\\ $\lb\left[\R_n(\psi_{k_1})\right]_{ij},\left[\R_n(\psi_{k_2})\right]_{i'j'}\rb^\hf{\T}$ forms a sesquilinear pair 
\[ \resizebox{\textwidth}{!}{
	$\lb U(\ell)=\frac{1}{\sqrt{n}}\lb \x(\ell)^{\T}  \A_n\y(\ell)-\rho(\ell)\tr(\A_n)\rb,~V(\ell')=\frac{1}{\sqrt{n}}\lb \x(\ell')^{\T}  \B_n\y(\ell')-\rho(\ell')\tr(\B_n)\rb\rb^{\T}$,}
\]
\noindent
 where
 \begin{gather*}
 A_n=\lb \psi_{k_1}\I_n-\frac{1}{n}\Y_2^{\T}  \V_{p'}\Y_2\rb^{-1}, ~B_n=\lb \psi_{k_2}\I_n-\frac{1}{n}\Y_2^{\T}  \V_{p'}\Y_2\rb^{-1},\\
 \rho(\ell)=\mathbb{E}(X_{\ell 1}Y_{\ell 1})=\delta_{ij},~\quad \rho(\ell')=\delta_{i'j'},
 \end{gather*}
 and $\x(\ell)$ corresponds to the $i$-th row of $\Y_1$, $\y(\ell)$ corresponds to the $j$-th row of $\Y_1$,
 $\x(\ell')$ corresponds to the $i'$-th row of $\Y_1$, $\y(\ell')$ corresponds to the $j'$-th row of $\Y_1$.	
 
Therefore
\[
\resizebox{\textwidth}{!}{
$w_3= \lim_{n\rightarrow \infty}\frac1n\tr\lb\lb \psi_{k_1}\I_n-\frac{1}{n}\Y_2^{\T}  \V_{p'}\Y_2\rb^{-1}\odot\lb \psi_{k_2}\I_n-\frac{1}{n}\Y_2^{\T}  \V_{p'}\Y_2\rb^{-1}\rb= \s(\psi_{k_1})\s(\psi_{k_2}),$}\]
\[
\resizebox{\textwidth}{!}{
$\theta_3= \lim_{n\rightarrow \infty}\frac1n\tr\lb\lb \psi_{k_1}\I_n-\frac{1}{n}\Y_2^{\T}  \V_{p'ii}\Y_2\rb^{-1}\lb \psi_{k_2}\I_n-\frac{1}{n}\Y_2^{\T}  \V_{p'}\Y_2\rb^{-1}\rb=\frac{ \s(\psi_{k_1})-\s(\psi_{k_2})}{\psi_{k_1}-\psi_{k_2}}.$}
\]
By Lemma \ref{lem:jointRn}, for any $1\leq i,j,i',j' \leq m$,
$\lb\left[\R_n(\psi_{k_1})\right]_{ij},\left[\R_n(\psi_{k_2})\right]_{i'j'}\rb^\hf{\T}$ weakly converges to a zero-mean Gaussian vector $\lb\left[\R(\psi_{k_1})\right]_{ij},\left[\R(\psi_{k_2})\right]_{i'j'}\rb^\hf{\T}$ with the following covariance structure:
\begin{itemize}
	\item[(1)] for $i=i',~j=j',~i\neq j$, $a_1=a_2=1, ~a_3=0$,
	\[\cov\lb\left[\R(\psi_{k_1})\right]_{ij},\left[\R(\psi_{k_2})\right]_{ij}\rb=\theta_3=\frac{ \s(\psi_{k_1})-\s(\psi_{k_2})}{\psi_{k_1}-\psi_{k_2}};\]
	\item[(2)] for $i=i'=j=j'$, $a_1=\nu_4-1,~a_2=1, ~a_3=1$,
	\[
	\resizebox{\textwidth}{!}{$
	\cov\lb\left[\R(\psi_{k_1})\right]_{ii},\left[\R(\psi_{k_2})\right]_{ii}\rb=2\theta_3+(\nu_4-3)w_3=\frac{2\lb \s(\psi_{k_1})-\s(\psi_{k_2})\rb}{\psi_{k_1}-\psi_{k_2}}+\beta_y\s(\psi_{k_1})\s(\psi_{k_2});$}\]
	\item[(3)] for all the other cases, $a_1=a_2=a_3=0$,
	\[\cov\lb\left[\R(\psi_{k_1})\right]_{ij},\left[\R(\psi_{k_2})\right]_{i'j'}\rb=0.\]
\end{itemize}
Then substituting $\s(\psi_{k})$ with \eqref{eq:s}, the limiting distribution of \\
$\left(
\sqrt{n}\lb \frac{\lambda_1}{\psi_1}-1\rb, \cdots, \sqrt{n}\lb \frac{\lambda_m}{\psi_m}-1\right)
\rb^\hf{\T}$  in Theorem \ref{thm:extreme} naturally follows.

\medskip\noindent{\bf \textsc{Step 2: Marginal limiting distribution of $\tr(\S_{22})$ and $\tr(\S_{11})$.}}
	
\noindent	
In this step we study the marginal limiting distribution of $\tr(\S_{22})$.  In fact, 
\[
\tr(\S_{22})=\tr\lb \frac{1}{n}\Y_2^{\T}\V_{p'}\Y_2\rb=\frac{1}{n}\sum_{i=1}^n \y_i^{\T}\V_{p'}\y_i,
\]
where each $\y_i$ is a random vector with $p'$ i.i.d. entries $Y_{ij}\lb 1\leq j\leq p'\rb$ satisfying $\E Y_{ij}=0$, $\E Y_{ij}^2=1$, $\E Y_{ij}^4=\nu_4$.

Moreover, by some calculations, we have
\[
\E\lb \tr(\S_{22})\rb=\frac{1}{n}\sum_{i=1}^n\E\lb \y_i^{\T}\V_{p'}\y_i\rb=\tr\lb \V_{p'}\rb,
\]
\begin{align*}
\E\lb \tr(\S_{22})\rb^2&=\frac{1}{n^2}\sum_{i\neq j}^n \E\lb \y_i^{\T}\V_{p'}\y_i\y_j^{\T}\V_{p'}\y_j\rb+\frac{1}{n^2}\sum_{i=1}^n\E\lb\y_i^{\T}\V_{p'}\y_i\rb^2\\
&=\frac{n^2-n}{n^2}\tr^2(\V_{p'})+\frac{1}{n}\left[\tr^2(\V_{p'})+2\tr(\V_{p'}^2)+(\nu_4-3)\sum_{i=1}^{p'}\left[\V_{p'}\right]_{ii}^2\right].
\end{align*}
Thus 
\[\var\lb \tr(\S_{22})\rb=\frac{2}{n}\tr(\V_{p'}^2)+\frac{\nu_4-3}{n}\sum_{i=1}^{p'}\left[\V_{p'}\right]_{ii}^2.\]
Actually, in Section 4.2.3 of \cite{BB16}, the authors have proved the asymptotic normality for trace of any symmetric polynomial of a general class of sample (auto-)covariance matrices. It is directly applicable to our case of $\tr(\S_{22})$ since our model settings fulfill all their assumptions. Therefore, we have, as $p, ~n\rightarrow \infty$, $p/n\rightarrow y>0$,	
\[\tr(\S_{22})-\tr(\V_{p'})\xrightarrow{d}\EN(0, ~2y\gamma_2+y(\nu_4-3)\gamma_{d,2}),\]
where $\gamma_2=\lim_{p'\rightarrow \infty}\frac{1}{p'}\tr(\V_{p'}^2),~\gamma_{d,2}=\lim_{p'\rightarrow \infty}\frac{1}{p'}\sum_{i=1}^{p'}[\V_{p'}]_{ii}^2$. 

Similarly, we have
\begin{equation*}
\E\lb \tr(\S_{11})\rb=\tr(\mathbf{\Lambda}),~\var\lb \tr(\S_{11})\rb=\frac{2}{n}\tr(\mathbf{\Lambda}^2)+\frac{\nu_4-3}{n}\sum_{i=1}^m\left[\mathbf{\Lambda}\right]_{ii}^2.
\end{equation*}
By Linderberg-Feller Central Limit Theorem, as $n\rightarrow \infty$, 
\begin{equation}\label{eq:S11}
\sqrt{n}\lb\tr(\S_{11})-\tr(\mathbf{\Lambda})\rb\xrightarrow{d}\EN\lb 0, ~ \sum_{i=1}^m\Lambda_{ii}^2(\nu_4-1)+\sum_{i\neq j} \Lambda_{ij}^2\rb.
\end{equation}

%
%

\noindent{\bf\textsc{Step 3: Joint limiting distribution of $\lb \lb\lambda_k\rb_{1\leq k\leq m},	\tr(\S_{n})\rb^{\T}$.}}
	
\noindent	
First,  by \eqref{eq:S11} , we have
\[\resizebox{\textwidth}{!}{$
\tr(\S_n)-\tr(\mathbf{\Sigma}_p)=\tr(\S_{11})-\tr\mathbf{ \Lambda}+\tr(\S_{22})-\tr(\V_{p'})=\tr(\S_{22})-\tr(\V_{p'})+O_{\mathbb{P}}\lb\frac{1}{\sqrt{n}}\rb.$}\]
Thus $\tr(\S_n)-\tr(\mathbf{\Sigma}_p)$ shares the same Gaussian limiting distribution with $\tr(\S_{22})-\tr(\V_{p'})$, i.e., under Assumptions (i)
 	to (iv),
	\[\tr(\S_n)-\tr(\mathbf{\Sigma}_p)\xrightarrow{d}\EN(0, ~2y\gamma_2+y(\nu_4-3)\gamma_{d,2}).\]

Secondly, from the previous proof we know that the main fluctuation of $\lambda_k's$ originates from $\R_n(\ell)$. It can be seen that $\lb\lambda_k\rb_{1\leq k\leq m}$  are asymptotically independent of $\Y_2$ because $\R_n(\ell)$ weakly converges \hf{to} a Gaussian matrix with distribution independent of $\Y_2$. 
 Actually,  by going through the proof of Theorem 11.10 (see Lemma \ref{prop:Rn} in this paper) and the result of Theorem 10.8 in \citet{YZB15a},  it can be proved that, conditioning on $\mathbf{Y}_2$, 
the limiting distribution of $\R_n(\ell)$ is a function of the LSD of $\S_{22}$, which does not depend on the value of the conditioning variable $\Y_2$. 

This establishes the asymptotic independence between $\R_n(\ell)$ and $\Y_2$. Moreover, since
{\small
	\[\R_n(\ell)=\frac{1}{\sqrt{n}}\lb \Y_1\lb \ell\I_n-\frac{1}{n}\Y_2^{\T}  \V_{p'}\Y_2\rb^{-1}\Y_1^{\T}  -\I_m\tr\lb \ell\I_n-\frac{1}{n}\Y_2^{\T}  \V_{p'}\Y_2\rb^{-1}\rb,
	\]}
if we treat $\frac{1}{\sqrt{n}}\Y_1\lb \ell\I_n-\frac{1}{n}\Y_2^{\T}  \V_{p'}\Y_2\rb^{-1}\Y_1^{\T} $ as $f(\Y_1,\Y_2)$, then
\begin{gather*}
\frac{1}{\sqrt{n}}\I_m\tr\lb \ell\I_n-\frac{1}{n}\Y_2^{\T}  \V_{p'}\Y_2\rb^{-1}=\E\lb f(\Y_1,\Y_2)|\Y_2\rb,\\
\R_n(\ell)=f(\Y_1,\Y_2)-\E\lb f(\Y_1,\Y_2)|\Y_2\rb, ~\cov(\R_n(\ell),\Y_2)=\mathbf{0}.
\end{gather*}

Note that $\tr(\S_{22})=\frac{1}{n}\tr\lb\Y_2^{\T}{\V_{p'}}\Y_2\rb$, the randomness of $\tr(\S_{22})$ all originates from $\Y_2$ and marginally $\tr(\S_{22})$ is also asymptotically normal.
Accordingly, we have $\lb \lb\lambda_k\rb_{1\leq k\leq m},\tr(\S_{22})\rb^\hf{\T}$ are asymptotically normal and independent.

Since $\tr(\S_{11})$ is of order $O(1/\sqrt{n})$ and $\lb\lambda_k\rb_{1\leq k\leq m}$ is of constant order, 
$\lb \lb\lambda_k\rb_{1\leq k\leq m},	\tr(\S_{n})\rb^\hf{\T}$ is also asymptotically normally distributed and the covariance between $\lb\lambda_k\rb_{1\leq k\leq m}$ and $\tr(\S_n)$ is of order $O(1/\sqrt{n})$. 

\subsection{Asymptotic joint distribution of $\lb\lb\lambda_k\rb_{1\leq k\leq m},\tr(\S_{11})\rb^{\T}$}

In this section we present a result of the joint distribution of $\lb\lb\lambda_k\rb_{1\leq k\leq m},\tr(\S_{11})\rb^{\T}$. It is crucial for quantifying second order terms  in Theorem \ref{prop:joint}.

\hf{Without loss of generality}, we consider $\lb\lambda_k,\tr(\S_{11})\rb^{\T}$ first. Since 
\[\tr(\S_{11})-\tr({\bf \Lambda})=\frac{1}{n}\tr\lb \Y_1^{\T}  {\bf \Lambda}\Y_1\rb-\tr({\bf \Lambda}),\]
the leading term of $\delta_{n,k}=\sqrt{n}\lb\frac{\lambda_k}{\psi_k}-1\rb$ is
\begin{equation*}
\frac{\left[\U^{\T}  \R_n(\psi_k)\U\right]_{kk}}{\s'(\psi_k)\psi_k}+\frac{\sqrt{n}(\s(\psi_k)-\s_n(\psi_k))}{\psi_k\s'(\psi_k)}.
\end{equation*}
The second term $\s_n(\psi_k)$ is a function of $\Y_2$ which is independent from $\tr(\S_{11})$. Therefore, we only have to consider the correlation between $\R_n(\psi_k)$ and $\tr(\S_{11})$. 

Note that $\tr(\S_{11})=\frac1n\tr({\bf \Lambda}\Y_1\Y_1^{\T}  )$ can be seen as linear combinations of entries in the  $m\times m$ matrix $\frac{1}{n}\Y_1\Y_1^{\T}  $. According to Lemma \ref{lem:jointRn}, \\
$\lb \R_n(\ell),~\frac{1}{\sqrt{n}}\lb\Y_1\Y_1^{\T}  -\I_m\rb\rb^{\T}$ forms a random sesquilinear pair with 
$$\A_n=\lb \ell\I_n-\frac{1}{n}\Y_2^{\T}  \V_{p'}\Y_2\rb^{-1},~\B_n=\I_n.$$
\hf{If} the correlation between each entry of $\R_n(\ell)$ and $\frac{1}{\sqrt{n}}\lb\Y_1\Y_1^{\T}  -\I_m\rb$ can be obtained, then we can derive the joint distribution of $\lb\lambda_k,\tr(\S_{11})\rb^{\T}$. 
More specifically, we have the following result, \hf{whose proof is relegated to the supplementary file}.

\begin{proposition}\label{lem:S11joint}
	Under Assumptions (i)-(iv), as $p,n\rightarrow \infty$, $p/n\rightarrow y$, we have
	\[
	\renewcommand{\arraystretch}{1.5}
	\lb \begin{array}{c}
	\sqrt{n}\lb \frac{\lambda_k}{\psi_k}-1\rb\\
	\sqrt{n}\lb\tr(\S_{11})-\tr({\bf \Lambda})\rb
	\end{array}\rb\xrightarrow{d} \EN \lb \lb
	\begin{array}{c}
	0\\
	0
	\end{array}
	\rb,~ \lb
	\begin{array}{cc}
	\sigma_1^2 & \rho_k\\
	\rho_k & \sigma_2^2
	\end{array}
	\rb\rb,
	\]
	where  
	\begin{gather*}
	\sigma_1^2=\sum_{i=1}^mu_{ik}^4\sigma_{\alpha_k}^2+\sum_{i\neq j}^mu_{ik}^2u_{jk}^2s_{\alpha_k}^2,~ 
	\sigma_2^2=\sum_{i=1}^m\Lambda_{ii}^2(\nu_4-1)+\sum_{i\neq j} \Lambda_{ij}^2,\\
	\rho_k=\frac{\alpha_k\psi'(\alpha_{k})}{\psi_k}\lb (\nu_4-1)\sum_{i=1}^m\Lambda_{ii}u_{ik}^2+\sum_{i\neq j}^m\Lambda_{ij}u_{ik}u_{jk}\rb,
	\end{gather*}
	${\bf \Lambda}=(\Lambda_{ij})_{m\times m}$,  has spectral decomposition ${\bf \Lambda}=\U\diag(\alpha_1,\cdots,\alpha_m)\U^{\T}  $, $\U=(u_{ij})_{m\times m}$, 
	$\sigma_{\alpha_k}^2$ and $s_{\alpha_k}^2$ are defined in \eqref{eq:sigmak} and \eqref{eq:sk}.
\end{proposition}

\section*{Acknowledgement}

{The authors thank the co-editor, Richard Samworth, the associate editor, and two referees for many constructive and motivating comments which greatly improve the paper. The research of F. Han was supported in part by NSF grant DMS-1712536.  J. Yao thanks support from  Hong Kong RGC grant GRF-17306918.
}

\newpage{}
{\renewcommand{\thesection}{A\arabic{section}}
\renewcommand{\theequation}{A\arabic{equation}}
\setcounter{section}{0}

This supplement contains the proofs of Theorem \ref{prop:joint} and Proposition \ref{lem:S11joint}, and also extra simulation results in Section~\ref{sec:App}.

\section{Proof of Theorem \ref{prop:joint}}

\noindent
In this proof we give more detailed characterization of the joint distribution of $\lb\lb\lambda_k\rb_{1\leq k\leq m}, \tr(\S_n)\rb^{\T}$, specifically \hf{regarding} some second order terms which are ignored in the proof of Theorem \ref{thm:extreme}. These second order terms are sorted out step by step as follows.

\medskip\noindent{\bf{Second order term for $\lb\lambda_k\rb_{1\leq k\leq m}$}}

\hspace{0.3cm}

\noindent
By \eqref{eq:delta} in the proof of Theorem \ref{thm:extreme}, we have
\[
\delta_{n,k}=\sqrt{n}\left(\frac{\lambda_{k}}{\psi_{k}}-1\right)= \frac{\left[\U^{\T}\R_n(\psi_{k})\U\right]_{kk}}{\psi_{k}\s'(\psi_{k})}-\frac{M_n(\psi_{k})}{\sqrt{n}\psi_{k}\s'(\psi_{k})}.
\]

\begin{lemma}\label{prop:Mz}
	(\cite{Zheng15})
	Define $M_n(z)=n(\s_{n}(z)-\s(z))$. Under Assumptions (i) to (iv), $M_n(z)$ converges weakly to a two-dimensional Gaussian process $M(z)$ satisfying
	\[
	\E(M(z))=\cfrac{y\int\frac{\s(z)^{3}t^{2}}{\left(1+t\s(z)\right)^{3}}\d H(t)}{\left(1-y\int\frac{t^{2}\s(z)^{2}}{\left(1+t\s(z)\right)^{2}}\d H(t)\right)^{2}}+\cfrac{y\beta_{y}\int\frac{\s(z)^{3}t^{2}}{\left(1+t\s(z)\right)^{3}}\d H(t)}{1-y\int\frac{t^{2}\s(z)^{2}}{\left(1+t\s(z)\right)^{2}}\d H(t)},
	\]
	and variance-covariance function
	\begin{align*}
	\cov(M(z_{1}),M(z_{2}))~&=~2\left(\frac{\s'(z_{1})\s'(z_{2})}{\left(\s(z_{1})-\s(z_{2})\right)^{2}}-\frac{1}{\left(z_{1}-z_{2}\right)^{2}}\right)\\
	~&+~y\beta_{y}\s'(z_{1})\s'(z_{2})\int\frac{t}{\left(1+t\s(z_{1})\right)^{2}}\cdot\frac{t}{\left(1+t\s(z_{2})\right)^{2}}\d H(t).
	\end{align*}
\end{lemma}

\noindent
For two sequence $\xi_k,~\eta_k,~k=1,2,\cdots$, denote $\xi_k \asymp \eta_k$ if \hf{there} exist two constants $c_1\geq c_2>0$ \hf{such that}
$$c_2\leq \liminf_{k\rightarrow \infty} \left\lvert\frac{\xi_k}{\eta_k}\right\rvert\leq \limsup_{k\rightarrow \infty} \left\lvert\frac{\xi_k}{\eta_k}\right\rvert\leq c_1.$$
Then by Lemma \ref{prop:Mz}, we have
\begin{align*}
\E\left(\frac{M(\psi_k)}{\sqrt{n}\psi_k\s'(\psi_k)}\right)&=\frac{\s(\psi_k)^{3}}{\sqrt{n}\psi_k\s'(\psi_k)}\left[\cfrac{y\int\frac{t^{2}}{\left(1-t/\alpha_k\right)^{3}}\d H(t)}{\left(1-y\int\frac{t^{2}}{\left(t-\alpha_k\right)^{2}}\d H(t)\right)^{2}}+\cfrac{y\beta_{y}\int\frac{t^{2}}{\left(1-t/\alpha_k\right)^{3}}\d H(t)}{1-y\int\frac{t^{2}}{\left(t-\alpha_k\right)^{2}}\d H(t)}\right]\\
~&\asymp\frac{1}{\sqrt{n}\alpha_k\psi_k},
\end{align*}
and from the variance-covariance function, we have
\begin{align*}
\var(M(z_2))~&=~\lim_{z_1\rightarrow z_2}\cov(M(z_1),M(z_{2}))=\lim_{z_1\rightarrow z_2}2\left(\frac{\s'(z_{1})\s'(z_{2})}{\left(\s(z_{1})-\s(z_{2})\right)^{2}}-\frac{1}{\left(z_{1}-z_{2}\right)^{2}}\right)\\
~&+\lim_{z_1\rightarrow z_2} y\beta_{y}\s'(z_{1})\s'(z_{2})\int\frac{t}{\left(1+t\s(z_{1})\right)^{2}}\cdot\frac{t}{\left(1+t\s(z_{2})\right)^{2}}\d H(t).
\end{align*}
By Taylor expansion,
\begin{gather*}
\s(z_1)-\s(z_2)=\s'(z_2)(z_1-z_2)+\frac{1}{2}\s''(z_2)(z_1-z_2)^2+\frac{1}{6}\s'''(z_2)(z_1-z_2)^3+o(z_1-z_2)^3,\\
\s'(z_1)-\s'(z_2)=\s''(z_2)(z_1-z_2)+\frac{1}{2}\s'''(z_2)(z_1-z_2)^^2+o(z_1-z_2)^2.
\end{gather*}
Accordingly, we obtain
\begin{align*}
~&\lim_{z_1\rightarrow z_2}2\left(\frac{\s'(z_{1})\s'(z_{2})}{\left(\s(z_{1})-\s(z_{2})\right)^{2}}-\frac{1}{\left(z_{1}-z_{2}\right)^{2}}\right)=2\lim_{z_1\rightarrow z_2}\cfrac{\s'(z_{1})\s'(z_{2})\left(z_{1}-z_{2}\right)^{2}-\left(\s(z_{1})-\s(z_{2})\right)^{2}}{\left(\s(z_{1})-\s(z_{2})\right)^{2}\left(z_{1}-z_{2}\right)^{2}}\\
= &~ 2\lim_{z_1\rightarrow z_2}\cfrac{\s'(z_{1})\s'(z_{2})-\left[\s'(z_2)+\frac{1}{2}\s''(z_2)(z_1-z_2)+\frac{1}{6}\s'''(z_2)(z_1-z_2)^2+o(z_1-z_2)^2\right]^2}{\left(\s(z_{1})-\s(z_{2})\right)^{2}}\\
~&~\lb \mbox{use expansion of } \s(z_1)-\s(z_2)\rb\\
= &~2\lim_{z_1\rightarrow z_2} \cfrac{s'(z_2)(\s'(z_{1})-\s'(z_{2}))-\s'(z_{2})\s''(z_{2})(z_1-z_2)}{\left(\s(z_{1})-\s(z_{2})\right)^{2}}-\cfrac{3(\s''(z_2))^2+4\s'(z_2)\s'''(z_2)}{12(\s'(z_2))^2}\\
~&~\lb \mbox{use both expansions }\rb\\
=&~ \cfrac{2\s'(z_2)\s'''(z_2)-3(\s''(z_2))^2}{6(\s'(z_2))^2},
\end{align*}
and hence
\[\var(M(z))=\cfrac{2\s'(z)\s'''(z)-3(\s''(z))^2}{6(\s'(z))^2}+y\beta(\s'(z))^2\int\frac{t^2}{(1+t\s(z))^4}\d H(t).\]
\noindent
Thus
\[\var\left(\frac{M(\psi_k)}{\sqrt{n}\psi_k\s'(\psi_k)}\right)=\cfrac{2\s'(\psi_k)\s'''(\psi_k)-3(\s''(\psi_k))^2}{6n\psi_k^2(\s'(\psi_k))^4}+\frac{y\beta}{n\psi_k^2}\int\frac{t^2}{(1+t\s(\psi_k))^4}\d H(t).\]
By \eqref{eq:s}, the leading term of
$2\s'(\psi_k)\s'''(\psi_k)-3(\s''(\psi_k))^2$ is $-\cfrac{4\psi''(\alpha_k)\s''(\psi_k)}{\alpha_k^2(\psi'(\alpha_k))^3}$.

\noindent
In conclusion, the second order term for   $\sqrt{n}\lb\frac{\lambda_k}{\psi_k}-1\rb$ is the Gaussian variable $\frac{M(\psi_k)}{\sqrt{n}\psi_k\s'(\psi_k)}$ and 
\[\E\left(\frac{M(\psi_k)}{\sqrt{n}\psi_k\s'(\psi_k)}\right)\asymp\frac{1}{\sqrt{n}\alpha_k\psi_k},~\var\left(\frac{M(\psi_k)}{\sqrt{n}\psi_k\s'(\psi_k)}\right)\asymp \frac{1}{n\psi_k^2}.\]

\medskip\noindent{\bf{Second order term for $\tr(\S_n)$}}

\hspace{0.3cm}

\noindent
Notice that
\[\tr(\S_n)-\tr(\mathbf{\Sigma}_p)=\tr(\S_{11})-\tr\mathbf{ \Lambda}+\tr(\S_{22})-\tr(\V_{p'}),\]
\[
\tr(\S_{22})-\tr(\V_{p'})\xrightarrow{d}\EN\left(0,~2y\gamma_2+y(\nu_4-3)\gamma_{d,2}\right),
\]
\[\sqrt{n}\lb\tr(\S_{11})-\tr({\bf\Lambda})\rb\xrightarrow{d}\EN\lb0,~\sum_{i=1}^m\Lambda_{ii}^2(\nu_4-1)+\sum_{i\neq j} \Lambda_{ij}^2\rb,\]
then the second order term for $\var(\tr(\S_n))$ is $\frac1n\lb\sum_{i=1}^m\Lambda_{ii}^2(\nu_4-1)+\sum_{i\neq j} \Lambda_{ij}^2\rb$.

\medskip\noindent{\bf{Second order term for $\cov\lb\lambda_k,\tr(\S_n)\rb$}}

\hspace{0.3cm}

\noindent
From \eqref{eq:delta} in the proof of Theorem \ref{thm:extreme}, we know that the leading term of $\delta_{n,k}=\sqrt{n}\lb\frac{\lambda_k}{\psi_k}-1\rb$ is
\begin{equation}\label{eq:delta1}
\frac{[\U^{\T}  \R_n(\psi_k)\U]_{kk}}{\s'(\psi_k)\psi_k}+\frac{\sqrt{n}(\s(\psi_k)-\s_n(\psi_k))}{\psi_k\s'(\psi_k)},
\end{equation}
where
\[
	\R_n(\ell)=\frac{1}{\sqrt{n}}\lb\Y_1\lb \ell\I_n-\frac{1}{n}\Y_2^{\T}  \V_{p'}\Y_2\rb^{-1}\Y_1^{\T}  -\tr\lb \ell\I_n-\frac{1}{n}\Y_2^{\T}  \V_{p'}\Y_2\rb^{-1}\I_m\rb,\]
and
\begin{align}\label{eq:l2}
~&~\tr(\S_{22})-\tr(\V_{p'})=\int x ~d \left[ p\lb F^{\S_{22}}(x)-F^{y_n,H_n}(x)\rb\right]\\ \nonumber
~=&-\frac{1}{2\pi i}\oint_{\mathcal{C}} z\cdot n(\s_n(z)-\s(z))\d z.
\end{align}
Here the contour $\mathcal{C}$ is closed and taken in the positive direction in the complex plane, enclosing the support of $F^{y,H}$, LSD of $\S_{22}$. The covariance between  $\delta_{n,k}$ and $\tr(\S_{22})-\tr(\V_{p'})$ originates from two parts:
\[\cov\lb \frac{[\U^{\T}  \R_n(\psi_k)\U]_{kk}}{\s'(\psi_k)\psi_k}, ~-\frac{1}{2\pi i}\oint_{\mathcal{C}} z\cdot n(\s_n(z)-\s(z))\d z\rb\]
and
\[ \cov\lb\frac{\sqrt{n}(\s(\psi_k)-\s_n(\psi_k))}{\psi_k\s'(\psi_k)},~-\frac{1}{2\pi i}\oint_{\mathcal{C}} z\cdot n(\s_n(z)-\s(z))\d z\rb.\]

\noindent
As for the first part, since $\R_n(\psi_k)$ is an $m\times m$ matrix, for $1\leq u,v\leq m$, we focus on
\[\cov\lb \left[\R_n(\psi_k)\right]_{uv},~n\s_n(z)\rb.\]
 Denote $\H(\ell):=\lb \ell\I_n-\frac{1}{n}\Y_2^{\T}  \V_{p'}\Y_2\rb^{-1}$ and $\H_{ij}(\ell)$ to be the $(i,j)$-th entry of $\H(\ell)$. Then, for $u=v$,
\[\left[\R_n(\ell)\right]_{uu}=\frac{1}{\sqrt{n}}\lb\sum_{i=1}^n(\y_{ui}^2-1)\H_{ii}(\ell)+\sum_{i\neq j}^n\y_{ui}\y_{uj}\H_{ij}(\ell)\rb,~ n\s_n(z)=-\sum_{i=1}^n \H_{ii}(z).\]
Notice that $\Y_1=\lb \y_{uk}\rb_{1\leq u\leq m,~1\leq k\leq n}$ is independent of $\Y_2$. Thus
\[\cov \lb \left[\R_n(\psi_k)\right]_{uu},~n\s_n(z)\rb= \frac{1}{\sqrt{n}}\lb\sum_{i,i'=1}^n\E(\y_{ui}^2-1)\H_{ii}(\psi_k)\H_{i'i'}(z)+\sum_{i\neq j,k}^n\E \y_{ui}\y_{uj}\H_{ij}(\psi_k)\H_{kk}(z)\rb=0.\]
Similarly for $u\neq v$,
\begin{gather*}
\left[\R_n(\ell)\right]_{uv}=\frac{1}{\sqrt{n}}\sum_{i,j=1}^n\y_{ui}\y_{vj}\H_{ij}(\ell),\\
\cov \lb \left[\R_n(\psi_k)\right]_{uv},~n\s_n(z)\rb=\frac{1}{\sqrt{n}}\sum_{i,j,i'=1}^n\E \y_{ui}\y_{vj}\H_{ij}(\psi_k)\H_{i'i'}(z)=0.
\end{gather*}
Therefore
\[\cov\lb \frac{[\U^{\T}  \R_n(\psi_k)\U]_{kk}}{\s'(\psi_k)\psi_k}, ~-\frac{1}{2\pi i}\oint_{\mathcal{C}} z\cdot n(\s_n(z)-\s(z))\d z\rb=0.\]

As for the second part, we have
\begin{align*} ~&~\cov\lb\frac{\sqrt{n}(\s(\psi_k)-\s_n(\psi_k))}{\psi_k\s'(\psi_k)},~-\frac{1}{2\pi i}\oint_{\mathcal{C}} z\cdot(\s_n(z)-\s(z))\d z\rb\\
=~&~\lb \frac{1}{\sqrt{n}\psi_k\s'(\psi_k)}\rb\cdot\frac{1}{2\pi i}\oint_{\mathcal{C}} z\cdot \cov(M_n(z),M_n(\psi_k)) \d z.
\end{align*}
Then by Lemma \ref{prop:Mz}, the covariance between the limits of $\frac{\sqrt{n}(\s(\psi_k)-\s_n(\psi_k))}{\psi_k\s'(\psi_k)}$ and $-\frac{1}{2\pi i}\oint_{\mathcal{C}} z\cdot(\s_n(z)-\s(z))\d z$  is
\begin{align*}
~&~\lb \frac{1}{\sqrt{n}\psi_k\s'(\psi_k)}\rb\cdot\frac{1}{2\pi i}\oint_{\mathcal{C}} z\cdot \cov(M(z),M(\psi_k)) \d z\\
~=&~~\lb \frac{1}{\sqrt{n}\psi_k\s'(\psi_k)}\rb\cdot\frac{1}{2\pi i}\oint_{\mathcal{C}} z\cdot 2\left[\left(\frac{\s'(z)\s'(\psi_k)}{\left(\s(z)-\s(\psi_k)\right)^{2}}-\frac{1}{\left(z-\psi_k\right)^{2}}\right)\right.\\
~&+~\left.y\lb \nu_4-3\rb\s'(z)\s'(\psi_k)\int\frac{t}{\left(1+t\s(z)\right)^{2}}\cdot\frac{t}{\left(1+t\s(\psi_k)\right)^{2}}\d H(t) \right] \d z.
\end{align*}
\noindent
More specifically, choos\hf{ing the} contour $\mathcal{C}$ with two poles outside $\mathcal{C}$:  $\s=\s(\psi_k)=-\frac{1}{\alpha_k}$ and $\s=0$, we have
\[\frac{2}{2\pi i}\oint_{\mathcal{C}} z\cdot \left(\frac{\s'(z)\s'(\psi_k)}{\left(\s(z)-\s(\psi_k)\right)^{2}}-\frac{1}{\left(z-\psi_k\right)^{2}}\right)\d z=~2y\s'(\psi_k)\int\frac{t^2}{(1+t\s(\psi_k))^2}\d H(t),\]
and
\begin{align*}
~&~\frac{1}{2\pi i}\oint_{\mathcal{C}} z\cdot y\lb \nu_4-3\rb\s'(z)\s'(\psi_k)\lb \int\frac{t}{\left(1+t\s(z)\right)^{2}}\cdot\frac{t}{\left(1+t\s(\psi_k)\right)^{2}}\d H(t)\rb \d z\\
~=&~y\lb \nu_4-3\rb\s'(\psi_k)\int\left[\frac{1}{2\pi i}\oint_{\mathcal{C}_{k}}\cfrac{t^2\lb -\frac{1}{\s}+y\int\frac{t}{1+t\s}\d H(t)\rb}{\left(1+t\s\right)^{2}\left(1+t\s(\psi_k)\right)^{2}}\d \s\right] \d H(t)\\
~=&~y\lb \nu_4-3\rb\s'(\psi_k)\int\frac{t^2}{(1+t\s(\psi_k))^2}\d H(t).
\end{align*}
In conclusion,  the leading term of covariance between $\sqrt{n}\lb\frac{\lambda_k}{\psi_k}-1\rb$ and $	\tr(\S_{22})-\tr(\V_{p'})$ is 
\[\frac{y(\nu_4-1)}{\sqrt{n}\psi_k}\int\frac{t^2}{(1+t\s(\psi_k))^2}\d H(t).\]
Meanwhile, due to Proposition \ref{lem:S11joint}, the covariance between  $\sqrt{n}\lb\frac{\lambda_k}{\psi_k}-1\rb$ and $	\tr(\S_{11})-\tr(\mathbf{\Lambda})$ is
\[\frac{\alpha_k\psi'(\alpha_{k})}{\sqrt{n}\psi_k}\lb (\nu_4-1)\sum_{i=1}^m\Lambda_{ii}u_{ik}^2+\sum_{i\neq j}^m\Lambda_{ij}u_{ik}u_{jk}\rb.\]
Thus the leading term of covariance between limits of $\sqrt{n}\lb\frac{\lambda_k}{\psi_k}-1\rb$ and $	\tr(\S_{n})-\tr(\mathbf{\Sigma}_p)$ is 
\[\frac{\alpha_k\psi'(\alpha_{k})}{\psi_k\sqrt{n}}\lb (\nu_4-1)\sum_{i=1}^m\Lambda_{ii}u_{ik}^2+\sum_{i\neq j}^m\Lambda_{ij}u_{ik}u_{jk}\rb+\frac{y(\nu_4-1)}{\sqrt{n}\psi_k}\int\frac{t^2}{(1+t\s(\psi_k))^2}\d H(t).\]

\section{Proof of Proposition \ref{lem:S11joint}}

\begin{proof}
	Denoting $\S_1=\frac{1}{n}\Y_1\Y_1^{\T}$, $\Y_1=(Y_{ij})_{m\times n}$, then as $n\rightarrow \infty$, we have,
	\begin{gather*}
	\mbox{for } 1\leq i\leq m,~	\sqrt{n}\lb[\S_1]_{ii}-1\rb=\frac{1}{\sqrt{n}}\sum_{k=1}^n \lb y_{ik}^2-1\rb\xrightarrow{d} \EN(0,\nu_4-1),\\ \mbox{ and for } 1\leq i\neq j\leq m,~\sqrt{n}[\S_1]_{ij}=\frac{1}{\sqrt{n}}\sum_{k=1}^ny_{ik}y_{jk}\xrightarrow{d} \EN(0,1).
	\end{gather*}
	Thus we have,
	\begin{align}\label{eq:S11}
	~&~\sqrt{n}\lb\tr(\S_{11})-\tr({\bf\Lambda})\rb=\sqrt{n}\lb\sum_{i,j=1}^m\Lambda_{ij}[\S_1]_{ij}-\tr({\bf \Lambda})\rb\\ \nonumber
	~&~\xrightarrow{d}\EN\lb 0, ~\sum_{i=1}^m\Lambda_{ii}^2(\nu_4-1)+\sum_{i\neq j} \Lambda_{ij}^2\rb.
	\end{align}
	
	
	Secondly,  for any $1\leq i,j,i',j' \leq m$, consider a pair of the $(i,j)$-th entry of $\R_n(\psi_k)$ and the $(i',j')$-th entry of $\sqrt{n}\lb\frac{1}{n}\Y_1\Y_1^{\T}  -\I_m\rb$, i.e., $\lb\left[\R_n(\psi_{k})\right]_{ij},\left[\sqrt{n}\lb\frac{1}{n}\Y_1\Y_1^{\T}  -\I_m\rb\right]_{i'j'}\rb^{\T}$. According to Lemma \ref{lem:jointRn}, this pair forms a random sesquilinear pair 
	\[\resizebox{\textwidth}{!}{$
		\lb
		U(\ell)=\frac{1}{\sqrt{n}}\lb \x(\ell)^{\T}  \A_n\y(\ell)-\rho(\ell)\tr(\A_n)\rb,~V(\ell')=\frac{1}{\sqrt{n}}\lb \x(\ell')^{\T}  \B_n\y(\ell')-\rho(\ell')\tr(\B_n)\rb\rb^{\T},$}
	\]
	where
	\begin{gather*}
	\A_n=\lb \psi_k\I_n-\frac{1}{n}\Y_2^{\T}  \V_{p'}\Y_2\rb^{-1},~\B_n=\I_n,\\
	\rho(\ell)=\mathbb{E}(X_{\ell 1}Y_{\ell 1})=\delta_{ij},~\quad \rho(\ell')=\delta_{i'j'},
	\end{gather*}
	and $\x(\ell)$ corresponds to the $i$-th row of $\Y_1$, $\y(\ell)$ corresponds to the $j$-th row of $\Y_1$,
	$\x(\ell')$ corresponds to the $i'$-th row of $\Y_1$, $\y(\ell')$ corresponds to the $j'$-th row of $\Y_1$.	
	
	Therefore, 
	\begin{gather*}
	w_3= \theta_3= \lim_{n\rightarrow \infty}\frac1n\tr\lb\lb \psi_{k}\I_n-\frac{1}{n}\Y_2^{\T}  \V_{p'}\Y_2\rb^{-1}\rb=-\s(\psi_k),\\
	\cov\lb U(\ell),V(\ell')\rb= w_3a_1+(\theta_3-w_3)(a_2+a_3)=w_3a_1.
	\end{gather*}
	Since $\rho(\ell)=\mathbb{E}(X_{\ell 1}Y_{\ell 1}),~	a_1=\mathbb{E}(X_{\ell 1}Y_{\ell 1}X_{\ell'1}Y_{\ell'1})-\rho(\ell)\rho(\ell')$, by Lemma \ref{lem:jointRn}, 
	$$\lb\left[\R_n(\psi_{k})\right]_{ij},\left[\sqrt{n}\lb\frac{1}{n}\Y_1\Y_1^{\T}  -\I_m\rb\right]_{i'j'}\rb^{\T}$$ weakly converges to a zero-mean Gaussian vector $\lb\left[\R(\psi_{k})\right]_{ij},[{\bf S}_0]_{i'j'}\rb^{\T}$ with the following covariance structure:
	\begin{itemize}
		\item[(1)] for $i=i',~j=j',~i\neq j$, $a_1=1$, $\rho(\ell)=\rho(\ell')=0$,
		\[\cov\lb\left[\R(\psi_{k_1})\right]_{ij},[{\bf S}_0]_{ij}\rb=-\s(\psi_k);\]
		\item[(2)] for $i=i'=j=j'$, $a_1=\nu_4-1$, $\rho(\ell)=\rho(\ell')=1$,
		\[\cov\lb\left[\R(\psi_{k_1})\right]_{ii},[{\bf S}_0]_{ii}\rb=-(\nu_4-1)\s(\psi_k);\]
		\item[(3)] for all other cases, $a_1=0$,
		\[\cov\lb\left[\R(\psi_{k})\right]_{ij},[{\bf S}_0]_{i'j'}\rb=0.\]
	\end{itemize}
	Note that as $p,n\rightarrow \infty,~p/n\rightarrow y$, $\R_n(\psi_k)$ weakly converges to  $\R(\psi_k)$, $\sqrt{n}\lb\frac{1}{n}\Y_1\Y_1^{\T}  -\I_m\rb$ weakly converges to ${\bf S}_0$ where  both $\R(\psi_k)$ and ${\bf S}_0$ are $m\times m$ random matrices with independent Gaussian entries. \hf{Notice that} the leading term of $\delta_{n,k}=\sqrt{n}\lb\frac{\lambda_k}{\psi_k}-1\rb$ is $
	\frac{\left[\U^{\T}  \R_n(\psi_k)\U\right]_{kk}}{\s'(\psi_k)\psi_k}$ and 
	\[\sqrt{n}\lb\tr(\S_{11})-\tr({\bf\Lambda})\rb=\tr\lb\sqrt{n}\lb\frac{1}{n}\Y_1\Y_1^{\T}  -\I_m\rb{\bf \Lambda}\rb.\]
	Therefore $\lb\sqrt{n}\lb\frac{\lambda_k}{\psi_k}-1\rb,~ \sqrt{n}\lb\tr(\S_{11})-\tr({\bf\Lambda})\rb\rb^\hf{\T}$ weakly converges to a two-dimensional zero-mean Gaussian vector $\lb \frac{\left[\U^{\T}  \R(\psi_k)\U\right]_{kk}}{\s'(\psi_k)\psi_k},~\tr({\bf \Lambda S}_0)\rb^\hf{\T}$ with covariance 
	\begin{align*}
	~&~\cov\lb \frac{\left[\U^{\T}  \R(\psi_k)\U\right]_{kk}}{\s'(\psi_k)\psi_k},~\tr({\bf \Lambda S}_0)\rb\\
	= &~\sum_{i=1}^m \cov\lb \frac{ u_{ik}^2\left[\R(\psi_k)\right]_{ii}}{\s'(\psi_k)\psi_k},\Lambda_{ii}[{\bf \S}_0]_{ii}\rb
	+\sum_{i\neq j}^m\cov\lb\frac{u_{ik}u_{jk}\left[\R(\psi_k)\right]_{ij}}{\s'(\psi_k)\psi_k},\Lambda_{ij}[{\bf \S}_0]_{ij}\rb\\
	=&~-\frac{\s(\psi_k)(\nu_4-1)}{\s'(\psi_k)\psi_k}\sum_{i=1}^m\Lambda_{ii}u_{ik}^2-\frac{\s(\psi_k)}{\s'(\psi_k)\psi_k}\sum_{i\neq j}^m\Lambda_{ij}u_{ik}u_{jk}.
	\end{align*}
	Substituting \eqref{eq:s} for $\s(\psi_{k})$ and $\s'(\psi_k)$ and combining with  \eqref{eq:S11}, Proposition \ref{lem:S11joint} follows.
\end{proof}

\section{Simulation results in Section~\ref{sec:App}}

Empirical data are generated following the factor model in
\eqref{eq:factor}. We consider the case of two factors $m=2$ with the
noise variance  $\sigma^2=2$. The signal-to-noise ratios for the two
distant spikes are parametrized as $t_k=\frac{\alpha_k}{\sigma^2},
~k=1,2$, where $t_2\geq c$ corresponds to test size, and
$1\leq t_2<c$ shows the test power. Different $c$ values are chosen to
illustrate the performance of our testing procedure
\eqref{eq:rejectreal} for detecting factors with various strength
levels. For each $(p,n)$ combination, the smallest $c$ is set 
slightly larger than $1+\sqrt{p/n}$, the signal-to-noise ratio of
the weakest factor that can be detected. Empirical sizes and powers are
shown in Tables \ref{tab:size}  and \ref{tab:power} for various
$(p,n)$ configurations and alternatives based on 3,000
replications. The nominal level is $\alpha=0.05$. As seen from Tables
\ref{tab:size} and \ref{tab:power}, our test shows better
performance for higher signal strength level $c$.

To demonstrate the necessity of
incorporating the second order correction to the asymptotic variance
proposed  in \eqref{eq:factorcorr}, numerical  comparison is made
between the testing procedures \eqref{eq:rejectreal} and \eqref{eq:reject2}. 
 Tables~\ref{tab:corr} and \ref{tab:corr2} give empirical results of this comparison  in both the above two-factor model and a newly added four-factor model with the same setup. It is observed that the
 refined variance function leads to a superior procedure in term of
 detection power.

\renewcommand{\arraystretch}{1.5}

\begin{table}[h!]
\centering
	\caption{Empirical size of \eqref{eq:rejectreal} for testing $H_0:~t_2=\frac{\alpha_2}{\sigma^2} \geq c$ based on 3,000 replications. $t_2\geq c$ corresponds to test size.  }\label{tab:size}
	\resizebox{0.8\textwidth}{!}{
			\begin{tabular}{cccc|cccc}
				\hline 
				$(p,n)$ & $c$ & $(t_{1},t_{2})$ & Size (5\%) & $(p,n)$ & $c$ & $(t_{1},t_{2})$ & Size (5\%)\tabularnewline
				\hline 
				$(100,200)$	& $c=3.5$ & (10, 3.5) & 0.054 & $(200,400)$ & $c=3.5$ & (10, 3.5) & 0.052 \tabularnewline
				&& (10, 3.7) &0.019 &&& (10, 3.7) &0.011 \tabularnewline
				& $c=5$ & (10, 5) &0.053 &  & $c=5$ & (10, 5) & 0.054 \tabularnewline
				&& (10, 5.5) & 0.007&&& (10, 5.5) &0.002 \tabularnewline 
				\hline 
				$(200,200)$& $c=3.5$ & (10, 3.5) & 0.044 &$(400,400)$  & $c=3.5$ & (10, 3.5) & 0.057 \tabularnewline
				&& (10, 3.7) &0.015  &&& (10, 3.7) & 0.011\tabularnewline
				& $c=5$ & (10, 5) & 0.057&  & $c=5$ & (10, 5) &0.058 \tabularnewline
				&& (10, 5.5) &0.005 &&& (10, 5.5) & 0.001\tabularnewline
				\hline 
				$(200,100)$& $c=3.5$ & (10, 3.5) & 0.029 &$(400,200)$  & $c=3.5$ & (10, 3.5) & 0.033 \tabularnewline
				&& (10, 3.7) &0.017  &&& (10, 3.7) &0.014 \tabularnewline 
				& $c=5$ & (10, 5) & 0.056 &  & $c=5$ & (10, 5) &0.050 \tabularnewline
				&& (10, 5.5) & 0.014 &&& (10, 5.5) & 0.006\tabularnewline
				\hline 
			\end{tabular}}
		\end{table}

		\begin{table}[h!]
		\centering
			\caption{Empirical power  of \eqref{eq:rejectreal} for testing $H_0:~t_2=\frac{\alpha_2}{\sigma^2} \geq c$ based on 3,000 replications. $t_2<c$ corresponds to test power. }\label{tab:power}
			\resizebox{0.8\textwidth}{!}{
					\begin{tabular}{cccc|cccc}
						\hline 
						$(p,n)$ & $c$ & $(t_{1},t_{2})$ & Power (5\%) & $(p,n)$ & $c$ & $(t_{1},t_{2})$ & Power (5\%)\tabularnewline
						\hline 
						$(100,200)$ & $c=3.5$ & (10, 3) & 0.427& $(200,400)$ & $c=3.5$ & (10, 3) &0.668 \tabularnewline
						&  & (10, 2.5) &0.923 &  &  & (10, 2.5) & 0.998\tabularnewline
						& $c=5$ & (10, 4) & 0.691 &  & $c=5$ & (10, 4) & 0.919\tabularnewline
						&  & (10, 3) & 1.000 &  &  & (10, 3) & 1.000 \tabularnewline
						\hline 
						$(200,200)$	& $c=3.5$ & (10, 3) & 0.345 &$(400,400)$  & $c=3.5$ & (10, 3) & 0.605 \tabularnewline
						&  & (10, 2.5) & 0.884 &  &  & (10, 2.5) & 0.993 \tabularnewline
						& $c=5$ & (10, 4) & 0.682 &  & $c=5$ & (10, 4) &0.914 \tabularnewline
						&  & (10, 3) &1.000  &  &  & (10, 3) & 1.000 \tabularnewline
						\hline 
						$(200,100)$ & $c=3.5$ & (10, 2.5) & 0.223 & $(400,200)$ & $c=3.5$ & (10, 2.5) &0.563 \tabularnewline
						&  & (10, 1.5) &0.450 &  &  & (10, 1.5) &0.834 \tabularnewline
						& $c=5$ & (10, 4) & 0.381 &  & $c=5$ & (10, 4) & 0.630 \tabularnewline
						&  & (10, 3) & 0.905  &  &  & (10, 3) &0.999 \tabularnewline
						\hline 
					\end{tabular}}
				\end{table}

\renewcommand{\arraystretch}{1.5}

\begin{table}[h!]
\centering
	\caption{Empirical size and power comparison between two  testing procedures \eqref{eq:rejectreal} and  \eqref{eq:reject2} for   $H_0:~t_2=\frac{\alpha_2}{\sigma^2} > c=5$, $T_{m_0}(\sigma_{*,2}^2)$ represents procedure \eqref{eq:rejectreal} after second order correction, $T_{m_0}(\tilde{\sigma}_2^2)$ corresponds to \eqref{eq:reject2} before correction.}\label{tab:corr}
	\resizebox{0.8\textwidth}{!}{
			\begin{tabular}{cccc|cccc}
				\hline 
				$(p,n)$ & $(t_{1},t_{2})$ &  $T_{m_0}(\sigma_{*,2}^2)$ &   $T_{m_0}(\tilde{\sigma}_2^2)$ & $(p,n)$ & $(t_{1},t_{2})$ & $T_{m_0}(\sigma_{*,2}^2)$ &   $T_{m_0}(\tilde{\sigma}_2^2)$ \tabularnewline
				\hline 
				$(100,200)$ & (10, 5) &0.052& 0.042  & $(200,400)$ & (10, 5) & 0.047  & 0.038 \tabularnewline
				& (10, 4) &0.700  &0.649 &  & (10, 4) & 0.920 & 0.905 \tabularnewline
				\hline 
				$(200,200)$ & (10, 5) & 0.053& 0.038  & $(400,400)$ & (10, 5) & 0.058 &0.049 \tabularnewline
				&(10, 4) & 0.670& 0.612 &  & (10, 4) & 0.922& 0.903\tabularnewline
				\hline 
				$(200,100)$ & (10, 5) &  0.053&0.026  & $(400,200)$ & (10, 5) & 0.058 &0.039 \tabularnewline
				&	 (10, 4) & 0.401 & 0.275 &  & (10, 4) & 0.621&0.539\tabularnewline
				\hline 
			\end{tabular}}
		\end{table}

		\begin{table}[h!]
		\centering
			\caption{Empirical size and power comparison between two  testing procedures \eqref{eq:rejectreal} and  \eqref{eq:reject2} for   $H_0:~t_4=\frac{\alpha_4}{\sigma^2} > c=5$, $T_{m_0}(\sigma_{*,4}^2)$ represents procedure \eqref{eq:rejectreal} after second order correction, $T_{m_0}(\tilde{\sigma}_4^2)$ corresponds to \eqref{eq:reject2} before correction.}\label{tab:corr2}
			\resizebox{0.8\textwidth}{!}{
					\begin{tabular}{cccc|cccc}
						\hline 
						$(p,n)$ & $(t_{1},t_{2},t_3,t_4)$ &  $T_{m_0}(\sigma_{*,4}^2)$ &   $T_{m_0}(\tilde{\sigma}_4^2)$ & $(p,n)$ & $(t_{1},t_{2},t_3,t_4)$ & $T_{m_0}(\sigma_{*,4}^2)$ &   $T_{m_0}(\tilde{\sigma}_4^2)$ \tabularnewline
						\hline 
						$(100,200)$ & (20,15,10,5) &0.078& 0.048  & $(200,400)$ & (20,15,10,5) &  0.065& 0.046\tabularnewline
						& (20,15,10,4) &  0.750&  0.657&  & (20,15,10,4) & 0.950 &0.925 \tabularnewline
						\hline 
						$(200,200)$ &  (20,15,10,5) & 0.077 & 0.044  & $(400,400)$ & (20,15,10,5) & 0.066  & 0.046 \tabularnewline
						&(20,15,10,4)& 0.743 & 0.631 &  & (20,15,10,4) & 0.933 & 0.903 \tabularnewline
						\hline 
						$(200,100)$ & (20,15,10,5)  & 0.082 & 0.026 & $(400,200)$ &  (20,15,10,5)  & 0.074 & 0.034\tabularnewline
						&	 (20,15,10,4) & 0.486  & 0.248&  & (20,15,10,4) & 0.706 & 0.561 \tabularnewline
						\hline 
					\end{tabular}}
				\end{table}

}


\begin{thebibliography}{99}
	
	\bibitem[Bai and Ding(2012)]{BD12}
	\textsc{Bai, Z.,  Ding, X. }(2012). Estimation of spiked eigenvalues in spiked models. \em{Random Matrices: Theory and Applications}, 1(2), 1150011.
	
	
	\bibitem[Bai and Silverstein(2004)]{BS04}
	\textsc{Bai, Z.,  Silverstein, J.}(2004). CLT for linear spectral statistics of large-dimensional sample covariance matrices. {\em The Annals of Probability}, 32(1A), 553-605.
	
	
	
	\bibitem[{Bai and Yao(2008)}]{BY08}
	\textsc{Bai, Z.,  Yao, J. }(2008). Central limit theorems for eigenvalues in a spiked population model.  {\em Annales de l'Institut Henri Poincar\'e, Probabilit\'es et Statistiques}, 44(3), 447-474.
	
	\bibitem[{Bai and Yao(2012)}]{BY12}
	\textsc{Bai, Z.,  Yao, J.}(2012). On sample eigenvalues in a generalized spiked population model. {\em Journal of Multivariate Analysis}, 106, 167-177. 
	
	\bibitem[Baik and Silverstein(2006)]{BS06}	
	\textsc{Baik, J., Silverstein, J. W. }(2006). Eigenvalues of large sample covariance matrices of spiked population models. {\em Journal of Multivariate Analysis}, 97(6), 1382-1408.
	
	\bibitem[Baik et al.(2005)]{BBP05}
	\textsc{Baik, J., Arous, G., P\'{e}ch\'{e}, S. }(2005). Phase transition of the largest eigenvalue for nonnull complex sample covariance matrices. {\em The Annals of Probability}, 33(5), 1643-1697.
	
	\bibitem[Bhattacharjee and Bose(2016)]{BB16}
	\textsc{Bhattacharjee, M., and Bose, A. }(2016). Large sample behaviour of high dimensional autocovariance matrices. {\em The Annals of Statistics}, 44(2), 598-628.
	
	
	\bibitem[Bianchi et al.(2011)]{BDMN11}
	\textsc{Bianchi, P., Debbah, M., Maida, M.,  Najim, J.} (2011). Performance of statistical tests for single-source detection using random matrix theory. {\em IEEE Transactions on Information Theory}, 57(4), 2400-2419.
	
	
	
	

\bibitem[Chen and Pan(2015)]{CP15}
\textsc{Chen, B., Pan, G.}(2015)  CLT for linear spectral statistics of normalized sample covariance matrices with the dimension much larger than the sample size.
{\em Bernoulli}, 21(2), 1089-1133. 


\bibitem[Choi et al(2017)]{Choi17}
\textsc{Choi, Y., Taylor, J., and Tibshirani, R. }(2017). Selecting the number of principal components: Estimation of the true rank of a noisy matrix. {\em The Annals of Statistics}, 45(6), 2590-2617.

	\bibitem[Chow and Teugels(1978)]{CT78}
	\textsc{Chow, T.,  Teugels, J.}(1978). The sum and the maximum of iid random variables. {\em In Proceedings of the 2nd Prague Symposium on Asymptotic Statistics},  81-92.
	
	\bibitem[Davis et al(2016)]{Davis16}
	\textsc{Davis, R. A., Heiny, J., Mikosch, T., and Xie, X.} (2016). Extreme value analysis for the sample autocovariance matrices of heavy-tailed multivariate time series. {\em Extremes}, 19(3), 517-547.
	
	\bibitem[Deo(2016)]{Deo16}
	\textsc{Deo, R.}(2016). On the Tracy–Widom approximation of studentized extreme eigenvalues of Wishart matrices. {\em Journal of Multivariate Analysis}, 147, 265-272.
	

	
	\bibitem[Hsing(1995)]{Hsing95}
	\textsc{Hsing, T. }(1995). A note on the asymptotic independence of the sum and maximum of strongly mixing stationary random variables. {\em The Annals of Probability}, 23(2), 938-947.
	
	\bibitem[Johnson and Graybill(1972)]{JG72}
	\textsc{Johnson, D.,  Graybill, F.}(1972). An analysis of a two-way model with interaction and no replication. {\em Journal of the American Statistical Association}, 67(340), 862-868.
	
	\bibitem[Johnstone(2001)]{Johnstone01}
	\textsc{Johnstone, I.}(2001). On the distribution of the largest eigenvalue in principal components analysis. {\em The Annals of Statistics}, 29(2), 295-327.
	
	\bibitem[Knowles and Yin(2017)]{KY14}
	\textsc{Knowles, A.,  Yin, J. }(2017). Anisotropic local laws for random matrices. {\em Probability Theory and Related Fields}, 169(1), 257-352.
	
	\bibitem[Kritchman and Nadler(2008)]{KN08}
	\textsc{Kritchman, S.,  Nadler, B. }(2008). Determining the number of components in a factor model from limited noisy data. {\em Chemometrics and Intelligent Laboratory Systems}, 94(1), 19-32.
	
	\bibitem[Ma(2012)]{Ma12}
	\textsc{Ma, Z. }(2012). Accuracy of the Tracy-Widom limits for the extreme eigenvalues in white Wishart matrices. {\em Bernoulli}, 18(1), 322-359.
	
	
	
	\bibitem[Nadler(2011)]{Nadler11}
	\textsc{Nadler, B. }(2011). On the distribution of the ratio of the largest eigenvalue to the trace of a Wishart matrix. {\em Journal of Multivariate Analysis}, 102(2), 363-371.
	
	
	\bibitem[Onatski et al.(2013)]{Onatski13}
	\textsc{Onatski, A., Moreira, M. J., and Hallin, M. }(2013). Asymptotic power of sphericity tests for high-dimensional data. {\em The Annals of Statistics}, 41(3), 1204-1231.
	
	\bibitem[Paul(2007)]{Paul07}
	\textsc{Paul, D. }(2007). Asymptotics of sample eigenstructure for a large dimensional spiked covariance model. {\em Statistica Sinica}, 17(4), 1617-1642.
	
	\bibitem[Paul and Aue(2014)]{PA14}
	\textsc{Paul, D.,  Aue, A. }(2014). Random matrix theory in statistics: A review. {\em Journal of Statistical Planning and Inference}, 150, 1-29.
	
	\bibitem[Silverstein and Choi(1995)]{SilvChoi95}
	\textsc{Silverstein, J.,  Choi, S.}(1995). Analysis of the limiting spectral distribution of large dimensional random matrices. {\em Journal of Multivariate Analysis}, 54(2), 295-309.	
	
	\bibitem[Wang and Fan(2017)]{FW15}
	\textsc{Wang, W., Fan, J. }(2017). Asymptotics of empirical eigenstructure for high dimensional spiked covariance model. {\em The Annals of Statistics}, 45(3), 1342-1374.
	
	\bibitem[Wang et al.(2014a)]{WSY14}
	\textsc{Wang, Q., Silverstein, J. W., Yao, J.}(2014a). A note on the CLT of the LSS for sample covariance matrix from a spiked population model. {\em Journal of Multivariate Analysis}, 130, 194-207.
	
	\bibitem[Wang et al.(2014b)]{WZY14}
	\textsc{Wang, Q., Su, Z., Yao, J.}(2014b). Joint CLT for several random sesquilinear forms with applications to large-dimensional spiked population models. {\em Electronic Journal of Probability}, 19(103), 1-28.
	
	\bibitem[Yao et al.(2015)]{YZB15a}
	\textsc{Yao, J., Zheng, S., Bai. Z.} (2015).
	Large Sample Covariance Matrices and High-Dimensional Data Analysis.
	{\em Cambridge University Press}.
	

	\bibitem[Zheng(2012)]{ZBY-cltF}
	\textsc{Zheng, S. }(2012). Central limit theorems for linear spectral statistics of large dimensional F-matrices. {\em Annales de l'Institut Henri Poincar\'e, Probabilit\'es et Statistiques}, 48(2), 444-476.
	
	
	\bibitem[Zheng et al.(2015)]{Zheng15}
	\textsc{Zheng, S., Bai, Z.,  Yao, J. }(2015). Substitution principle for CLT of linear spectral statistics of high-dimensional sample covariance matrices with applications to hypothesis testing. {\em The Annals of Statistics},  43(2), 546-591.
	


	
	
	 \bibitem[Zheng et al.(2016)]{ZBYZ16}
	\textsc{Zheng, S., Bai, Z., Yao, J., Zhu, H. }(2016). CLT for linear spectral statistics of large dimensional sample covariance matrices with dependent data. {\em  preprint} arXiv:1708.03749.
	
\end{thebibliography}
\end{document}